\documentclass[titlepage,11pt]{article}
\usepackage{latexsym}
\usepackage{amsfonts}
\usepackage{amsmath}
\usepackage{comment}
\usepackage{tikz}
\usepackage{float}
\newcommand\blackslug{\hbox{\hskip 1pt \vrule width 4pt height 8pt depth 1.5pt
        \hskip 1pt}}
\newcommand\bbox{\hfill \quad \blackslug \bigbreak}
\def\DD{\hbox{-}}
\def\CC{\hbox{-}\cdots\hbox{-}}
\def\LL{,\ldots,}
\def\cupcup{\cup\cdots\cup}
\newcommand{\vare}{\varepsilon}

\usepackage{url}

\title{Pure pairs. IX. Transversal trees}
\author{Alex Scott\thanks{Research supported by EPSRC grant
EP/V007327/1.}\\
Mathematical Institute, University of Oxford, Oxford OX2 6GG, UK
\\
\\
Paul Seymour\thanks{Supported by AFOSR grants FA9550-22-1-0234 and A9550-19-1-0187, and NSF
grants DMS-2154169 and DMS-1800053.}\\
Princeton University, Princeton, NJ 08544
\\
\\
Sophie Spirkl\thanks{We acknowledge the support of the Natural Sciences and Engineering Research
Council of Canada (NSERC), [funding reference number RGPIN-2020-03912]. Cette
recherche a \'et\'e financ\'ee par le Conseil de recherches en sciences naturelles et
en g\'enie du Canada (CRSNG), [num\'ero de r\'ef\'erence RGPIN-2020-03912].}\\
University of Waterloo, Waterloo, Ontario N2L3G1, Canada}

\date{September 1, 2020; revised \today\\ Accepted manuscript. The published version appeared in SIAM Journal on Discrete Mathematics, Volume 38(1), 2024, see \url{https://doi.org/10.1137/21M1456509}}

\newtheorem{thm}{}[section]

\newcommand{\Proof}{\noindent{\bf Proof.}\ \ }

\begin{document}
\maketitle
\begin{abstract}
Fix $k>0$, and let $G$ be a graph, with vertex set partitioned into $k$ subsets (``blocks'') of approximately equal size. An induced subgraph of $G$
is ``transversal'' (with respect to this partition) if it has exactly one vertex in each block (and therefore it has exactly $k$ vertices). 
A ``pure pair'' in $G$
is a pair $X,Y$ of disjoint subsets of $V(G)$ such that either all edges between $X,Y$ are present or none are; and in the present context we
are interested in pure pairs $(X,Y)$ where each of $X,Y$ is a subset of one of the blocks, and not the same block. This paper collects
several results and open questions concerning how large a pure pair must be present if various types of transversal subgraphs are excluded.

\end{abstract}
\section{Introduction}
Graphs in this paper are finite, and without loops or parallel edges. Let $A,B\subseteq V(G)$ be disjoint.
We say that $A$ is {\em complete} to $B$, or $A,B$ are {\em complete}, if every vertex in $A$ is adjacent to every vertex in $B$,
and similarly $A,B$ are {\em anticomplete} if no vertex in $A$ has a neighbour in $B$. 
A {\em pure pair} in $G$ is a pair $A,B$ of disjoint subsets of $V(G)$ such that $A,B$
are complete or anticomplete. The number of vertices of $G$ is denoted by $|G|$.
The complement graph of $G$ is denoted by $\overline{G}$.
If $X\subseteq V(G)$, $G[X]$ denotes the subgraph induced on $X$.

A {\em blockade} $\mathcal{B}$ in a graph $G$ is a sequence $(B_1\LL B_k)$ of pairwise disjoint nonempty subsets of $V(G)$, called the {\em blocks} of 
$\mathcal{B}$; and the {\em width} of $\mathcal{B}$ is the minimum cardinality of its blocks, and its {\em length} is $k$. An induced subgraph $H$ of $G$ is
{\em $\mathcal{B}$-rainbow} if each of its vertices belongs to a block of $\mathcal{B}$, and it has at most one vertex in each block;
and {\em $\mathcal{B}$-transversal} if it is $\mathcal{B}$-rainbow and has exactly one vertex in each block (and therefore
has exactly $k$ vertices). A {\em copy} of a graph $H$ in a graph $G$ means an induced subgraph of $G$ that is isomorphic to $H$.

In earlier papers of this sequence we proved several theorems that say that if we have a blockade $\mathcal{B}$ and there is no 
$\mathcal{B}$-rainbow copy of some special graph $H$, then there must be a pure pair $X,Y$ with $|X|,|Y|$ large (in terms of the width 
of the blockade). For instance, in~\cite{pure1} we proved:
\begin{thm}\label{pure1thm}
For every forest $H$, there exists $d>0$, such that,
for every graph $G$ with a blockade $\mathcal{B}$ of length at least $d$, if every vertex of $G$ has degree less than $W/d$, and there 
is no anticomplete pair $X,Y$ in $G$ with $|X|,|Y|\ge W/d$ (where $W$ is the width
of $\mathcal{B}$), then
there is a $\mathcal{B}$-rainbow copy of $H$ in $G$.
\end{thm}
In this paper we investigate what happens if we ask for a $\mathcal{B}$-transversal copy of $H$ rather than just a $\mathcal{B}$-rainbow copy.

This leads naturally to the question, what if we ask even more? Let $H$ be a graph together with a fixed linear ordering of 
its vertex set, say $v_1,v_2\LL v_k$. (We call this an ``ordered graph''.)
If $\mathcal{B}=(B_1\LL B_k)$ is a blockade in a graph $G$, an {\em ordered $\mathcal{B}$-transversal copy of $H$}
means a $\mathcal{B}$-transversal induced subgraph $J$ of $G$, such that there is an isomorphism $\phi$ from $H$ to $J$,
with $\phi(v_i)\in B_i$ for $1\le i\le k$.
Erd\H{o}s, Hajnal and Pach~\cite{EHP} proved:
\begin{thm}\label{EHPthm}
For every ordered graph $H$, there exist $c,\vare>0$, such that for every graph $G$, and every blockade $\mathcal{B}=(B_1\LL B_k)$ in $G$,
where $k=|H|$, either:
\begin{itemize}
\item there is an ordered $\mathcal{B}$-transversal copy of $H$ in $G$; or
\item there are distinct $i,j\in \{1\LL k\}$, and $X\subseteq B_i$ and $Y\subseteq B_j$, such that $X,Y$ is a pure pair, and 
$|X|, |Y|\ge \vare W^c$, where $W$ is the width of $\mathcal{B}$.
\end{itemize}
\end{thm}
This is where we will start. We would particularly like to know, for which graphs $H$ can we take  $c=1$ in the second 
bullet? (Sadly, almost never: only for graphs $H$ with at most two vertices.) And, more promising: 
for which $H$ can we get $|X|\ge \vare W$ and $|Y|\ge \vare W^c$ in the second bullet?

At the other extreme, we cannot get past the following open question, a variant of a conjecture of
Conlon, Fox and Sudakov~\cite{fox} that is discussed further in~\cite{pure3} (a ``triangle'' means a copy of the complete graph $K_3$):

\begin{thm}\label{foxconj}
{\bf Question: } Do there exist $c,\vare>0$ with the following property? Let $\mathcal{B}=(B_1,B_2,B_3)$ be a blockade in a graph $G$, 
and let $W$ be its width. If there is no $\mathcal{B}$-transversal triangle,
then there exist distinct $i,j\in \{1,2,3\}$ and $X\subseteq B_i$ and $Y\subseteq B_j$, such that $X,Y$ is a pure pair, and $|X|\ge\vare W$,
and $|Y|\ge \vare W^c$.
\end{thm}
Settling \ref{foxconj} was our initial goal in this research (although we were not able to do it), 
and it is easy to see that it can be reduced to the sparse case,
when for all 
distinct $i,j\in \{1,2,3\}$, every vertex in $B_i$ has at most $|B_j|/100$ neighbours in $B_j$ (and $100$ can be replaced by any other number).

We say a blockade $\mathcal{B}=(B_1\LL B_k)$ has {\em local degree} $\lambda$ if $\lambda\ge 0$ is the maximum of the number of neighbours of $v$ in $B_j$, maximized over 
all distinct $i,j\in \{1\LL k\}$ and all $v\in B_i$. (We set this to be zero if $k\le 1$.)
All our results concern blockades with local degree at most a small constant times the width. Also, we will look for pure pairs $X,Y$
where $X,Y$ are each a subset of a block, and not the same block; and we will no longer need pure pairs $X,Y$ with $X$ complete to $Y$.
Let us say a blockade $\mathcal{B}=(B_1\LL B_k)$ is
{\em $(x,y)$-cohesive} if for all distinct  $i,j\in \{1\LL k\}$, there do not exist $X\subseteq B_i$ and $Y\subseteq B_j$
such that $|X|\ge x$, and $|Y|\ge y$, and $X$ is anticomplete to $Y$. (This is true if $k\le 1$.)
Here are our main results, first for unordered copies of $H$:
\begin{thm}\label{cycle}
	(Proved in \ref{4cycle} and \ref{kcycle}.) If $H$ is a cycle with $k\ge 4$ vertices, then there exist $\vare,c>0$ with the following property. 
Let $\mathcal{B}=(B_1\LL B_k)$ be a 
blockade in a 
graph $G$, with local degree less than $\vare W$ and  $(\vare W,\vare W^c)$-cohesive, where $W$ is its width.
Then there is a $\mathcal{B}$-transversal copy of $H$ in $G$.
\end{thm}
This statement for $k=3$ is open and equivalent to \ref{foxconj}.

\begin{thm}\label{path}
	(Proved in \ref{betterpath}.) For every integer $k\ge 1$, there exists $\vare>0$ with the following property. Let $\mathcal{B}=(B_1\LL B_k)$ be a 
blockade in a 
graph $G$, with local degree less than $\vare W$ and $(\vare W,\vare W)$-cohesive, where $W$ is its width.
Then there is a $\mathcal{B}$-transversal copy of a $k$-vertex path in $G$.
\end{thm}
The statement of \ref{path} is also true for the tree obtained from a path with $k-2$ vertices by adding two extra 
vertices, both adjacent to the last vertex of the path,
provided that $k\ge 5$; and also for the tree obtained from a path with $k-t$ vertices by adding $t$ extra vertices, each 
adjacent to the last vertex of the
path, provided that $k\ge t^22^t$. (These are proved in \ref{brooms}.) It is {\em not} true for the tree obtained from a $(k-6)$-vertex path by 
adding six extra vertices, three adjacent to the first vertex of the path and 
three adjacent to the last. (This is \ref{doublebroom}.)

For ordered copies of $H$, we have:
\begin{thm}\label{tree}
	(Proved in \ref{tree2}.) If $H$ is an ordered tree with $k\ge 2$ vertices, then there exists $\vare>0$ with the following property. 
Let $\mathcal{B}=(B_1\LL B_k)$ be a 
blockade in a 
graph $G$, with local degree less than $\vare W$ and $(\vare W,\vare W^{1/(k-1)})$-cohesive, where $W$ is its width.
Then there is an ordered $\mathcal{B}$-transversal copy of $H$.
\end{thm}

For caterpillars we can strengthen this. A {\em caterpillar} is a tree in which all the vertices of degree more than one belong to one path.
\begin{thm}\label{caterpillar}
	(Proved in \ref{orderedcaterpillar2}.) If $H$ is an ordered caterpillar with $k$ vertices, then there exists $\vare>0$ with the following property. 
Let $\mathcal{B}=(B_1\LL B_k)$ be a       
blockade in a 
graph $G$, with local degree less than $\vare W$, and $(\vare W,\vare W^{1/d})$-cohesive, where $W$ is its width, and $d$ is the maximum degree of $H$.
Then there is an ordered $\mathcal{B}$-transversal copy of $H$.

\end{thm}

And a counterexample (if $t\ge 1$ is an integer, $S_t$ denotes the star with $t+1$ vertices, that is, the tree in which 
one vertex is adjacent to all the others):
\begin{thm}\label{counterexample}
	(Proved in \ref{orderedstar}.) Let $t\ge 3$ be an integer, and let $S_t^+$ be obtained from $S_t$ by linearly 
ordering its vertex set.
For all $c>1/t$ and all $\vare>0$, there is a graph $G$ and a 
blockade $\mathcal{B}=(B_1\LL B_{t+1})$
in $G$, with local degree less than $\vare W$ and $(\vare W,\vare W^c)$-cohesive where $W$ is its width,
such that there is no ordered $\mathcal{B}$-transversal copy of $S_t^+$ in $G$.
\end{thm}
\section{Two easy covering theorems}

We say that a graph $H$ has the {\em strong transversal property} if there exists $\vare>0$ with the following property:
for every graph $G$, if $\mathcal{B}=(B_1\LL B_{|H|})$ is a
blockade in 
$G$, with local degree less than $\vare W$ and $(\vare W,\vare W)$-cohesive, where $W$ is its width,
then there is a $\mathcal{B}$-transversal copy of $H$ in $G$. If this holds we say that $\vare$ is an {\em STP-coefficient}.
We start with asking, which graphs have the strong transversal property?
A sparse random graph of girth at least $|H|+1$ shows that every such graph $H$ must be a forest, and one might hope
that all forests have the property, extending the results of~\cite{pure1}, but this is not true, as we shall see.
Nevertheless, some forests have the strong transversal property: here is what we know about them.
\begin{thm}\label{strongbonsai}
Let $H$ be a graph. 
\begin{itemize}
\item If $H$ is not a forest then $H$ does not have the strong transversal property.
\item If every component of $H$ has the strong transversal property then so does $H$ (the converse is false).
\item If $H$ is a path,  then $H$ has the strong transversal property. 
\item If $|H|> 4$ and $H$ is obtained from a path by adding two new vertices adjacent to the last 
vertex of the path, then $H$ has the strong transversal property. 
\item If $|H|> 2^t(t^2-t+1)$ and $H$ is obtained from a path by adding $t$ new vertices adjacent to the last      
vertex of the path, then $H$ has the strong transversal property.
\item If $H$ is obtained from a path by adding six new vertices, three adjacent to the first vertex of the path and
three adjacent to the last, then $H$ does not have the strong transversal property.
\item If $H$ has a vertex of degree at least $d$ where $2^{d-1}\ge |H|$, then $H$ does not have the strong transversal property.
\end{itemize}
\end{thm}

In particular, it is {\em not} true that if a graph $H$ has the property then so do all its induced subgraphs, or 
indeed all its components:  the graph
$S_3$ does not have the property (by the last bullet of \ref{strongbonsai}) but if we add a vertex of degree zero,
this five-vertex forest has the property.
Indeed, it follows from one of the results of~\cite{pure1} that for any forest, if we add enough vertices of degree zero
we will obtain a forest with the strong transversal property. 

We will prove the various statements of \ref{strongbonsai}
as separate theorems (except for the first two, which we leave to the reader). 

If $G$ is a graph and $A,B\subseteq V(G)$ are disjoint, we say $A$ {\em covers} $B$ if every vertex in $B$ has a neighbour in $A$.
For convenience, let us say a blockade $\mathcal{B}=(B_1\LL B_k)$ in $G$ is {\em $\vare$-coherent} if for all distinct $i,j\in \{1\LL k\}$:
\begin{itemize}
\item each vertex in $B_i$ has fewer than $\vare |B_j|$ neighbours in $B_j$, and
\item there do not exist $X\subseteq B_i$ and $Y\subseteq B_j$ with $|X|\ge \vare |B_i|$
and $|Y|\ge \vare |B_j|$ and $X$ anticomplete to $Y$.
\end{itemize}

This is very much like $\vare$-cohesion, but is different when the blocks have different sizes. 
This new definition is not really needed, but it works nicely and is a little more compact than using cohesion 
(and we used the same concept in earlier papers).
Let us first prove the third statement of \ref{strongbonsai}, that is, \ref{path},
 which we restate in a slightly stronger form:

\begin{thm}\label{betterpath}
	Let $k\ge 2$ be an integer, and $0<\vare\le 1/(2k-2)$. Let $\mathcal{B}=(B_1\LL B_k)$ be an $\vare$-coherent blockade in a
graph $G$.
Then there is a $\mathcal{B}$-transversal $k$-vertex path in $G$ with an end-vertex in $B_1$.
\end{thm}
\Proof  We define $t_1\LL t_k$ with $\{t_1\LL t_k\}=\{1\LL k\}$, and $A_i\subseteq B_i$ for $1\le i\le k$, as follows.
Let $t_1=1$. 
Inductively, let $1\le i\le k$, and suppose that $t_1\LL t_{i}$ and $A_{t_1}\LL A_{t_{i-1}}$ have been defined, with the properties that
\begin{itemize}
	\item $\emptyset\ne A_{t_h}\subseteq B_{t_h}$ for $1\le h< i$;
	\item for $1\le h< i-1$, $A_{t_h}$ covers $A_{t_{h+1}}$;
	\item for $1\le g<h<  i$ with $h- g\ge 2$, there are no edges between $A_{t_g}$ and $A_{t_h}$;
	\item for each $j\in \{1\LL k\}\setminus \{t_1\LL t_{i}\}$, at least $(1-2(i-1)\vare)|B_j|$ vertices in $B_j$ have no neighbour
		in $A_{t_1}\cupcup A_{t_{i-1}}$; and
		\item if $i\ge 2$, at least $\vare|B_{t_{i}}|$ vertices in $B_{t_{i}}$ have a neighbour in $A_{t_{i-1}}$ 
			and have no neighbour in
                $A_{t_1}\cupcup A_{t_{i-2}}$.
\end{itemize}
Let $J=\{1\LL k\}\setminus \{t_1\LL t_{i}\}$.
For each $j\in J$, let $C_j$ be the set of vertices in $B_j$
with no neighbour in $A_{t_1}\cupcup A_{t_{i-1}}$; thus $|C_j|\ge (1-2(i-1)\vare)|B_j|$.
If $i = 1$ let $D=B_{t_{i}}$, and if $i\ge 2$ let $D$ be the set of vertices in $B_{t_{i}}$ that have a neighbour in 
$A_{t_{i-1}}$ and have no neighbour in
		$A_{t_1}\cupcup A_{t_{i-2}}$; thus $|D|\ge  \vare|B_{t_{i}}|$.  
If $i=k$, let $A_{t_{i}}=D$ and the inductive definition is complete, so we assume that $i\le k-1$.
For each $j\in J$, fewer than $\vare|B_j|$ vertices in $B_j$ have no neighbour in $D$, since $|D|\ge  \vare|B_{t_{i+1}}|$ and
the blockade is $\vare$-coherent; and since $|C_j|\ge (1-2(i-1)\vare)|B_j|\ge 2\vare|B_j|$, at least $\vare|B_j|$ vertices in $C_j$
have a neighbour in $D$.

Since $J\ne \emptyset$, there exists $A_{t_{i}}\subseteq D$ minimal such that for 
some $j\in J$, at least $\vare|B_{j}|$ vertices in $C_{j}$ have a neighbour in $A_{t_{i}}$. From the minimality
of $A_{t_{i}}$, for each $j\in J$ there are fewer than $2\vare|B_j|$ vertices in $C_j$ with a neighbour in $A_{t_{i}}$,
and hence there are at least $|C_j|-2\vare|B_j|\ge (1-2i\vare)|B_j|$ vertices in $C_j$ with no neighbour in $A_{t_{i}}$.
Choose $j\in J$ such that at least $\vare|B_{j}|$ vertices in $C_{j}$ have a neighbour in $A_{t_{i}}$, and define $t_{i+1}=j$.
This completes the inductive definition.

Choose $a_{t_k}\in A_{t_k}$. Since $A_{t_{k-1}}$ covers $A_{t_k}$, there exists
$a_{t_{k-1}}\in A_{t_{k-1}}$ adjacent to $a_{t_k}$; and similarly for $i=k-2,k-3\LL 1$ there exists $a_{t_i}\in A_{t_i}$ 
adjacent to $a_{t_{i+1}}$.
But for $1\le i<j\le k$ with $j\ge i+2$, there are no edges between $A_{t_i}, A_{t_j}$; so $a_{t_1}\DD a_{t_2}\CC a_{t_k}$ is an induced path.
This proves \ref{betterpath}.~\bbox

A somewhat similar proof (the proofs of \ref{betterpath} and \ref{findstar} are both specializations of the proof of the main theorem 
of~\cite{cats}) shows:
\begin{thm}\label{findstar}
Let $k\ge 1$ be an integer, let $K=2^{k-1}+1$, and let $0<\vare\le 3^{-K}$. Let $\mathcal{B}=(B_1\LL B_K)$ be an $\vare$-coherent
blockade in a
graph $G$.
Then there is a $\mathcal{B}$-rainbow copy of $S_k$ in $G$.
\end{thm}
\Proof
If $J\subseteq \{1\LL K\}$, a {\em star-partition} of $\mathcal{B}$ (see figure~\ref{fig:star-partition}) with {\em element set} $J$ is a sequence
$$(h_1,I_1,h_2,I_2\LL h_t, I_t)$$
where $J= \{h_1\LL h_t\}\cup I_1\cupcup I_t$, 
and a choice of a nonempty subset $A_i\subseteq B_i$ for each $i\in J$,
with the following properties:
\begin{itemize}
\item $t\ge 1$, and $h_1\LL h_t\in \{1\LL K\}$ are distinct, and the sets $I_1\LL I_t$ are pairwise disjoint subsets of 
$\{1\LL K\}\setminus \{h_1,h_2\LL h_t\}$;
\item for $1\le s\le t$, and all $j\in I_s$,  $A_j$ covers $A_{h_s}$;
\item for $1\le s\le t$ and all $j\in I_s$, $A_j$ is anticomplete to 
\begin{itemize}
\item all the sets $A_{j'}$ for $j'\in I_s\setminus \{j\}$;
\item all the sets $A_{h_{s'}}$ for $s'\in \{1\LL t\}\setminus \{s\}$; and
\item all the sets $A_{j'}$ for $s'\in \{1\LL t\}\setminus \{s\}$ and $j'\in I_{s'}$.
\end{itemize}
\end{itemize}
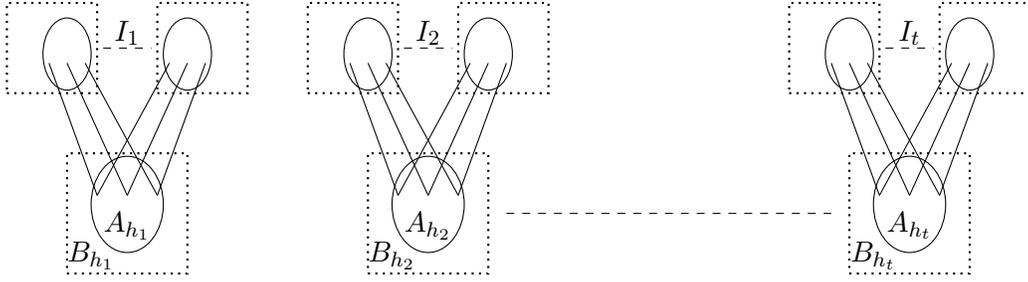
\begin{figure}[H]
	\centering
\begin{tikzpicture}[scale=0.8,auto=left]
	\draw[thick, dotted] (-8,-1) rectangle (-6,1);
\draw[thick, dotted] (-3,-1) rectangle (-1,1);
\draw[thick, dotted] (5,-1) rectangle (7,1);
	\draw[thick, dotted] (-9,2) rectangle (-7.5,3.5);
	\draw[thick, dotted] (-6.5,2) rectangle (-5,3.5);
	\draw[thick, dotted] (-4,2) rectangle (-2.5,3.5);
        \draw[thick, dotted] (-1.5,2) rectangle (0,3.5);
	\draw[thick, dotted] (4,2) rectangle (5.5,3.5);
        \draw[thick, dotted] (6.5,2) rectangle (8,3.5);
	\draw[thin, dashed] (-.7,0)--(4.7,0);
	\draw[thin, dashed] (-7.4,2.75)--(-6.6,2.75);
	\draw[thin, dashed] (-2.4,2.75)--(-1.6,2.75);
	\draw[thin, dashed] (5.6,2.75)--(6.4,2.75);
	\draw (-8,2.65) ellipse (.4 and .6);
	\draw (-6,2.65) ellipse (.4 and .6);
	\draw (-3,2.65) ellipse (.4 and .6);
	\draw (-1,2.65) ellipse (.4 and .6);
	\draw (5,2.65) ellipse (.4 and .6);
	\draw (7,2.65) ellipse (.4 and .6);
	\draw (-7,.15) ellipse (.6 and .8);
	\draw (-2,.15) ellipse (.6 and .8);
	\draw (6,.15) ellipse (.6 and .8);

	\draw (-8,2.5) -- (-7,.3);
	\draw (-8.3,2.5) -- (-7.5,.3);
	\draw (-7.7,2.5) -- (-6.5,.3);

	\draw (-6,2.5) -- (-7,.3);
        \draw (-6.3,2.5) -- (-7.5,.3);
        \draw (-5.7,2.5) -- (-6.5,.3);

	\draw (-3,2.5) -- (-2,.3);
        \draw (-3.3,2.5) -- (-2.5,.3);
        \draw (-2.7,2.5) -- (-1.5,.3);

	 \draw (-1,2.5) -- (-2,.3);
        \draw (-1.3,2.5) -- (-2.5,.3);
        \draw (-0.7,2.5) -- (-1.5,.3);

	\draw (5,2.5) -- (6,.3);
        \draw (4.7,2.5) -- (5.5,.3);
        \draw (5.3,2.5) -- (6.5,.3);

	\draw (7,2.5) -- (6,.3);
        \draw (6.7,2.5) -- (5.5,.3);
        \draw (7.3,2.5) -- (6.5,.3);

	\tikzstyle{every node}=[]
\node at (-7,-.2) {$A_{h_1}$};
\node at (-2,-.2) {$A_{h_2}$};
\node at (6,-.2) {$A_{h_t}$};
\node at (-7.6,-.7) {$B_{h_1}$};
\node at (-2.6,-.7) {$B_{h_2}$};
\node at (5.4,-.7) {$B_{h_t}$};
\node at (-7, 3) {$I_1$};
\node at (-2, 3) {$I_2$};
\node at (6, 3) {$I_t$};

\end{tikzpicture}
\caption{A star-partition.} \label{fig:star-partition}
\end{figure}

The {\em linkage} of a star-partition is the maximum over all distinct $i,j\in \{h_1\LL h_t\}$
of the maximum over $v\in A_i$ of $n/|B_j|$ where $n$ is the number of neighbours of $v$ in $A_j$ (or zero if $t\le 1$).
Its {\em length} is $t$, 
and its {\em value} is
$2^{|I_1|}+\cdots+2^{|I_t|}$.
There is a star-partition with linkage less than $\vare$, length $K$ and value $K$, since we may set $t=K$ and the 
sets $I_1\LL I_K$
all empty, and $A_i=B_i$ for $1\le i\le K$. 
Choose $t\ge 1$ minimum such that there is a star-partition with 
linkage less than $\vare3^{K-t}$, length $t$ and value at least $K$, say 
$$(h_1,I_1,h_2,I_2\LL h_t, I_t)$$
and $A_i\subseteq B_i$ for each $i\in J$, where $J$ is its element set.

Suppose that $t\ge 2$. We may assume that $|I_t|\le |I_1|\LL |I_{t-1}|$. Choose $C_{h_t}\subseteq A_{h_t}$
minimal such that $C_{h_t}$ covers at least one-third of one of the sets $A_{h_1}\LL A_{h_{t-1}}$, say of $A_{h_1}$.
Let $C_{h_1}$ be the set of vertices in $A_{h_1}$ with a neighbour in $C_{h_t}$, and for $2\le s\le t-1$, 
let $C_{h_s}$ be the set of vertices in $A_{h_t}$ with no neighbour in $C_{h_t}$. Thus $|C_{h_1}|\ge |A_{h_1}|/3$, and 
from the minimality of $C_{h_t}$,
and since the linkage is less than $\vare3^{K-t}$, it follows that
$$|C_{h_s}|\ge (2/3-\vare 3^{K-t})|A_{h_s}|\ge |A_{h_s}|/3$$
for $2\le s\le t-1$.
Let $J'=J\setminus I_t$, and for each $i\in J'$ with $i\notin \{h_1\LL h_{t-1}, h_t\}$ let $C_i=A_i$.
Then 
$$(h_1, I_1\cup \{h_t\}, h_2, I_2\LL h_{t-1}, I_{t-1})$$
and the sets $C_i\;(i\in J')$ form a star-partition with
linkage less than $\vare3^{K-t+1}$, length $t-1$ and value at least $K$, contrary to the minimality of $t$.

Hence $t=1$, and so $2^{|I_1|}\ge K$, and therefore $|I_1|\ge k$. We may assume that $h_1=1$, and $2\LL k+1\in I_1$.
Choose $u\in A_1$; then for $2\le i\le k+1$ there exists $v_i\in A_i$ adjacent to $u$, and the subgraph induced on
$\{u,v_1\LL v_k\}$ is a $\mathcal{B}$-rainbow copy of $S_k$. This proves \ref{findstar}.~\bbox

Next we show that the expression $2^{k-1}+1$ in \ref{findstar} is best possible, and hence the final statement of 
\ref{strongbonsai} holds, because of the following:
\begin{thm}\label{noclaw}
For every integer $k\ge 2$, and all $\vare>0$,
there is a graph $G$, and an $\vare$-coherent blockade $\mathcal{B}$ in $G$ of length $2^{k-1}$,
such that
there is no $\mathcal{B}$-rainbow copy of $S_k$ in $G$.
\end{thm}
The proof needs the following
two lemmas, which will also be needed later in the paper. The first is a standard estimate:
\begin{thm}\label{binom}
If $n\ge k\ge 1$ are integers then $\binom{n}{k}\le (en/k)^k$ (where $e$ is Euler's number).
\end{thm}
\Proof
By Stirling's formula~\cite{robbins}, we have 
$$k! \ge  \sqrt{2\pi} k^{k+1/2}e^{-k + 1/(12k+1)} \ge (k/e)^k$$ 
for $k\ge 1$, and so 
$$\binom{n}{k}\le n^k/k!\le n^k(k/e)^{-k}=(en/k)^k.$$
This proves \ref{binom}.~\bbox

The second lemma is also a well-known result.
\begin{thm}\label{random}
Let $\vare>0$; then there exists $d>0$ such that for all sufficiently large $n$,
there is a 
bipartite graph with bipartition $(A,B)$, where $|A|=|B|=n$, 
such that every vertex has degree less than $d$, and there do not
exist anticomplete sets $A'\subseteq A$ and $B'\subseteq B$ with $|A'|,|B'|\ge \vare n$.
\end{thm}
\Proof
Choose $c >4\vare^{-2}\log 2 $, and choose $d$ with  $2\log 2 < d\log(d/(2ce))$.
Now let $n$ be sufficiently large, and let $A,B$ be disjoint sets both of cardinality $2n$.
For each $u\in A$ and $v\in B$, let $u,v$ be adjacent with probability $c/n$, independently.
The expected number of anticomplete pairs $(A',B')$ with $A'\subseteq A$ and $B'\subseteq B$ and with 
$|A'|,|B'|\ge \vare n$ is at most 
$$2^{4n}(1-c/n)^{\vare^2n^2}\le 2^{4n}e^{-c\vare^2 n}\le 1/2$$
(since $c\vare^2>4\log 2$) if $n$ is sufficiently large. The probability that a given vertex $v\in A\cup B$ has degree 
at least $d$ is at most 
$$\binom{2n}{d}(c/n)^d\le (2en/d)^d(c/n)^d =(2ce/d)^d,$$ 
by \ref{binom}.
Consequently the probability that at least $n$ vertices in $A$ have degree
at least $d$ is at most $2^{2n}(2ce/d)^{nd}<1/4$ (since $2\log 2 < d\log(d/(2ce))$) if $n$ is sufficiently large. Hence, if
$n$ is sufficiently large, with
positive probability there is no anticomplete pair $(A',B')$ with $A'\subseteq A$ and $B'\subseteq B$ and with 
$|A'|,|B'|\ge \vare n$, and at most $n$ vertices in $A$, and at most $n$ vertices in $B$ have degree at least $d$.
Thus by deleting the $n$ vertices in $A$ with largest degree, and the same for $B$, we obtain a graph satisfying the 
theorem. This proves \ref{random}.~\bbox

To prove \ref{noclaw}, for inductive purposes we will prove something a little stronger, the following.
\begin{thm}\label{betterclaw}
Let $\vare>0$. For every integer $k\ge 2$, and every integer $p\ge 0$, there exists $W(k,p)$ such that for all 
integers $W\ge W(k,p)$,
if $H$ is a graph with 
$2^{k-1}W$ vertices, and with maximum degree at most $p$, there is a graph $G$ with the same vertex set and with $H$
as a subgraph, and an $\vare$-coherent
blockade $\mathcal{B}$ in $G$ of length $2^{k-1}$, such that
there is no $\mathcal{B}$-rainbow copy of $S_k$ in $G$.
\end{thm}
\Proof 
We prove \ref{betterclaw} by induction on $k$.
Suppose first that $k=2$. Choose $d$ and $n_0$ such that \ref{random} holds for all $n\ge n_0$.
Let $W(2,p)= \max(n_0, (p+d)/\vare)$; we claim that $W(2,p)$ satisfies the theorem.
Let $W\ge W(2,p)$, and let $H$ be a graph
 with
$2W$ vertices, and with maximum degree at most $p$. Let $B_1,B_2$ be two disjoint subsets of $V(H)$ both of 
cardinality $W$. By \ref{random} there is a graph $J$ with bipartition $(B_1,B_2)$, such that every vertex has degree
less than 
$d$, and there is no anticomplete pair $(A',B')$ with $A'\subseteq B_1$ and $B'\subseteq B_2$ and with
$|A'|,|B'|\ge \vare W$, that is, in $J$ the blockade  $\mathcal{B}=(B_1,B_2)$ has local degree less than $d$ and is 
$(\vare W,\vare W)$-cohesive. Let $G$ be the union of $H,J$; then in $G$ the same blockade $\mathcal{B}$
has local degree less than $p+d$ and is $(\vare W,\vare W)$-cohesive. Since it only has two blocks and therefore
there is no $\mathcal{B}$-rainbow copy of $S_2$, the result holds.

Now we assume inductively that $k\ge 3$ and the theorem holds for $k-1$. 
Choose $d$ and $n_0$ such that \ref{random} holds for all $n\ge n_0$, with $\vare$ replaced by $2^{2-k}\vare$.

Let 
$$W(k,p)=\max \left(n_0, (p+d)/\vare, W(k-1,p+(p+d)^2)\right).$$
We claim that $W(k,p)$ satisfies the theorem. Let $W\ge W(k,p)$, and
let $H$ be a graph
 with
$2^{k-1}W$ vertices, and with maximum degree at most $p$. Let $V_1,V_2$ be two disjoint subsets of $V(H)$ both of 
cardinality $2^{k-2}W$. By \ref{random} there is a graph $J$ with bipartition $(V_1,V_2)$, such that every vertex 
has degree less than 
$d$, and there is no anticomplete pair $(A',B')$ with $A'\subseteq V_1$ and $B'\subseteq V_2$ and with
$$|A'|,|B'|\ge (2^{2-k}\vare)|A'|=\vare W.$$
For $i = 1,2$, let $H_i$ be the graph with vertex set $V_i$, in which distinct vertices
$u,v\in V_i$ are adjacent if and only if either they are $H$-adjacent, or they have a common 
$(H\cup J)$-neighbour in $V(H)\setminus V_i$.
It follows that $H_i$ has maximum degree at most $p+(p+d)^2$. Since $W(k,p)\ge W(k-1,p+(p+d)^2)$, the inductive hypothesis implies
that 
for $i = 1,2$ there is a graph $G_i$ with vertex set $V_i$
and with $H_i$
as a subgraph, and an $\vare$-coherent blockade $\mathcal{B}_i$  
in $G_i$ of length $2^{k-2}$, such that
there is no $\mathcal{B}_i$-rainbow copy of $S_{k-1}$ in $G_i$. Let $G$ be the union of $G_1,G_2$ and $H\cup J$, and let 
$\mathcal{B}$ be the blockade with blocks all the blocks of $\mathcal{B}_1$ and all those of $\mathcal{B}_2$ 
(in some order). It follows that, in $G$, the blockade $\mathcal{B}$ is $\vare$-coherent (since 
$p+d\le \vare W$).
Suppose there is a $\mathcal{B}$-rainbow copy $T$ of $S_k$ in $G$, and let $v$ be the vertex
of $T_k$ that has degree $k$ in $T_k$. From the symmetry we may assume that $v\in V_1$. Since all neighbours 
of $v$ in $V_2$ are 
pairwise adjacent (since $H_2\subseteq G_2\subseteq G$) and $T$ is induced in $G$, it follows that 
at most one vertex of $T$
belongs to $V_2$, and so there is a $\mathcal{B}_1$-rainbow copy of $S_{k-1}$ in $G_1$, a contradiction. This proves \ref{betterclaw}.~\bbox

\section{Brooms}
Let $P$ be a path with vertices $p_1\LL p_k$ in order; and let $H$ be obtained from $P$ by adding $t$ new vertices, 
each adjacent to $p_k$. We define $B(k,t)=H$; such a graph is a {\em broom}.

If instead we add $s+t$ new vertices to $P$, $s$ of them adjacent to $p_1$ and the other $t$ to $p_k$, the graph we produce is 
called a {\em double broom} and is denoted $B(k,s,t)$. We still have three parts of \ref{strongbonsai} to prove (namely that $B(k,2)$
has the strong transversal property, and so does $B(k,t)$ if $k\gg t$, and that $B(k,3,3)$ does not), and we will do that 
in this section.
We begin with the easiest:
\begin{thm}\label{doublebroom}
For every integer $k\ge 1$, the double broom $B(k,3,3)$ does not have the strong transversal property.
\end{thm}
\Proof
Suppose that $B(k,3,3)$ has the strong transversal property, with STP-coefficient $\vare$.
Choose $d$ and $n_0$ such that \ref{random} holds for all $n\ge n_0$.
Let 
$$W=\max\left(n_0, \left(((k+3)d+2d+3(k+3)d^2)^4 + d +3(k+3)d^2\right)/\vare\right).$$
Let $B_1\LL B_{k+6}$ be pairwise disjoint sets each of cardinality $W$.
Let $\mathcal{B}=(B_1\LL B_{k+6})$, and $V_1=B_{k+4}\cup B_{k+5}\cup B_{k+6}$, and $V_2=B_1\cupcup B_{k+3}$.
For $1\le i<j\le k+6$, let $J_{i,j}$ be a graph with bipartition $(B_i, B_j)$ with maximum degree less than $d$, 
such that there is no anticomplete pair $A',B'$ with $A'\subseteq B_i$ and $B'\subseteq B_j$ and $|A'|,|B'|\ge \vare W$.
Let $J$ be the union of all the graphs $J_{i,j}$. Let $L$ be the graph with vertex
set $V_1$ in which distinct $u,v$ are adjacent if there is a $\mathcal{B}$-rainbow path in $J$ with ends $u,v$,
of length one or two and with its interior vertex (if any) in $V_2$. 
Let $R$ be the graph with vertex set $V_2$
in which distinct $u,v$ are adjacent if 
there is a $\mathcal{B}$-rainbow path of $J\cup L$ with ends $u,v$ and with every internal 
vertex in $V_1$. Let $G=J\cup L\cup R$.

Since $J$ is a subgraph of $G$ it follows that $\mathcal{B}$ is $(\vare W,\vare W)$-cohesive in $G$. The only edges of $G$ 
between $V_1,V_2$ are those of $J$; and $L$ has maximum degree at most $2d+3(k+3)d^2$, since each vertex in $V_1$
has degree at most $(k+3)d$ in $J$, and each of those neighbours has degree at most $3d$ in $J$. If $P$ is a 
$\mathcal{B}$-rainbow path of $J\cup L$ with ends in $V_2$ and with every internal 
vertex in $V_1$, then $P$ has at most five vertices; and since $J\cup L$ has maximum degree at most $(k+3)d+2d+3(k+3)d^2$, it
follows that each vertex in $V_2$ is an end of at most $((k+3)d+2d+3(k+3)d^2)^4$ such paths, and so $R$
has maximum degree at most  $((k+3)d+2d+3(k+3)d^2)^4$.

Consequently $\mathcal{B}$ has local degree at most 
$$(k+3)d+2d+3(k+3)d^2)^4 + d +3(k+3)d^2<\vare W.$$
From the choice of $\vare$, there is a $\mathcal{B}$-rainbow copy $H$ of $B(k,3,3)$ in $G$. Let it be constructed from a $k$-vertex path
with vertices $p_1\LL p_k$ in order, by adding three new vertices $q_1,q_2,q_3$ adjacent to $p_1$ and three new vertices 
$r_1,r_2,r_3$ adjacent to $p_k$. Every path of $H$ with both ends in $V_2$ has all its internal vertices in $V_2$, since
otherwise there would be an induced $\mathcal{B}$-rainbow path of $G$ with both ends in $V_2$ and all internal vertices in
$V_1$, with its ends nonadjacent; and this would contradict the construction of $R$. If one of $p_1\LL p_k$ belongs to $V_1$
then one of $q_1,q_2,q_3$ and one of $r_1,r_2,r_3$ belongs to $V_2$, contradicting the claim just made; so $p_1\LL p_k\in V_2$.
So exactly three of $q_1,q_2,q_3,r_1,r_2,r_3$ belong to $V_1$; and so two of them have a common neighbour in $V_2$,
contradicting the definition of $L$. This proves that there is no such $\vare$, and so proves \ref{doublebroom}.~\bbox

Let $\mathcal{B}=(B_1\LL B_k)$ be a blockade in $G$, and let $J$ be a digraph with vertex set $\{1\LL k\}$. 
(Digraphs in this paper do not have loops or parallel edges, but they may have antiparallel edges. Thus, if there is an edge 
from $i$ to $j$ then it is unique, 
but there might also be an edge from $j$ to $i$.)
For $\tau>0$, we say that
$J$ is a {\em $\tau$-covering digraph} for $\mathcal{B}$ if for $1\le i\le k$ there exists $A_i\subseteq B_i$, and for each edge $ij$
of $J$ there exists $X_{ij}\subseteq B_i$, with 
the following properties:
\begin{itemize}
\item $|A_i|\ge \tau^{|E(J)|}|B_i|$ for $1\le i\le k$; 
\item for each edge $ij$ of $J$,  $X_{ij}$ covers $A_j$ and $X_{ij}$ is anticomplete to $A_h$ for
all $h\in \{1\LL k\}\setminus\{i,j\}$; and
\item for all edges $ij, i'j'$ of $J$ with $i\ne i'$, and $i\ne j'$ and $i'\ne j$, the sets $X_{ij}, X_{i'j'}$ are anticomplete.
\end{itemize}
We call the sets $(A_1\LL A_k)$ a {\em core} for $J$.
There is a $\tau$-covering digraph, because we can take $J$ with no edges and $A_i=B_i$ for each $i$.
A $\tau$-covering digraph for $\mathcal{B}$ is {\em optimal} if no $\tau$-covering digraph for $\mathcal{B}$ has strictly more edges. 
If $X\subseteq V(G)$ we denote by $N(X)$ or $N_G(X)$ the set of vertices in $V(G)\setminus X$ that have a neighbour in $X$.

\begin{thm}\label{coveringdig}
Let $\mathcal{B}=(B_1\LL B_k)$ be a blockade in a graph $G$, let $\tau>0$, and let 
$J$ be an optimal $\tau$-covering digraph for $\mathcal{B}$, with core $(A_1\LL A_k)$. Suppose that
for all distinct $i,j\in \{1\LL k\}$, every vertex in $B_i$ has fewer than $(1-2\tau)\tau^{|E(J)|}|B_j|$ neighbours in $B_j$. 
Let $1\le i\le k$ and 
$X\subseteq A_i$. Then either 
\begin{itemize}
\item there exists $j\in \{1\LL k\}\setminus \{i\}$ such that $ij\in E(J)$ and $|A_j\setminus N(X)|<\tau^{|E(J)|+1}|B_j|$; or 
\item for every $j\in \{1\LL k\}\setminus \{i\}$ such that $ij\notin E(J)$ we have $|A_j\cap N(X)|< \tau^{|E(J)|+1}|B_j|$.
\end{itemize}
\end{thm}
\Proof
Let $s=|E(J)|$.
We assume the second bullet of the theorem is false, so there exists $j\in \{1\LL k\}\setminus \{i\}$ with $ij\notin E(J)$ 
such that $|A_j\cap N(X)|\ge \tau^{s+1}|B_j|$.
Choose $Y\subseteq X$ minimal such that there exists $j\in \{1\LL k\}\setminus \{i\}$ 
with $ij\notin E(J)$ such that $|A_j\cap N(Y)|\ge \tau^{s+1}|B_j|$.
Let $C_j=A_j\cap N(Y)$, let $C_i=A_i$,
and for each $h\in \{1\LL k\}\setminus \{i,j\}$ let $C_h=A_h\setminus N(Y)$. Thus
$|C_j|\ge \tau^{s+1}|B_j|$. If $|C_h|\ge \tau^{s+1}|B_h|$ for every $h\in \{1\LL k\}\setminus \{i,j\}$, then 
adding the edge $ij$ to $J$ gives a $\tau$-covering digraph
for $\mathcal{B}$ with core $(C_1\LL C_k)$, contrary to the optimality of $J$. Thus there exists $h\in \{1\LL k\}\setminus \{i,j]$
with $|C_h|<\tau^{s+1}|B_h|$. If $ih\notin E(J)$, the minimality of $Y$ 
and the hypothesis about local degree imply that 
$$|A_h\setminus C_h|<\left(\tau^{s+1}+\left((1-2\tau)\tau^s\right)\right)|B_h|=\left(\tau^s-\tau^{s+1}\right)|B_h|,$$
and so 
$$|C_h|\ge \tau^{s+1}|B_h|,$$
a contradiction. Thus $ih\in E(J)$. Since $|C_h|<\tau^{s+1}|B_h|$, and 
$$A_h\setminus N(X)\subseteq A_h\setminus N(Y)=C_h,$$
the first bullet of the theorem holds. 
This proves \ref{coveringdig}.~\bbox

We deduce:
\begin{thm}\label{morecoveringdig}
Let $\mathcal{B}=(B_1\LL B_k)$ be a blockade of length at least two in a graph $G$, and let $J$ be an optimal $\tau$-covering
digraph for $\mathcal{B}$, with core $(A_1\LL A_k)$. If $\mathcal{B}$ is $(1-2\tau)\tau^{|E(J)|}$-coherent,
then every vertex of $J$ has outdegree at least one.
\end{thm}
\Proof
Suppose that $i$ has outdegree zero in $J$; we may assume that $i=1$. Let $|E(J)|=s$. Since $k\ge 2$ and $\mathcal{B}$ is 
$(1-2\tau)\tau^{s}$-coherent and $|A_1|\ge \tau^{s}|B_1|$, it follows that $|A_2\setminus N(A_1)|<(1-2\tau)\tau^{s}|B_2|$,
and so 
$$|A_2\cap N(A_1)|> |A_2|-(1-2\tau)\tau^{s}|B_2|\ge \tau^{s+1}|B_2|.$$ 
But then both the outcomes of \ref{coveringdig} (with $i = 1$ and $X=A_1$) are false, a contradiction. This proves \ref{morecoveringdig}.~\bbox

We use these results to prove one of the remaining parts of \ref{strongbonsai}:

\begin{thm}\label{brooms}
Let $k,t\ge 0$ be integers with $t\ge 2$ and $k\ge 2^t(t^2-t+1)-t+1$. Then $B(k,t)$ has the strong transversal property.
\end{thm}
%
%
%
%
%
%
%
%
%
%
\Proof
Let $\tau=1/6$, and 
let $\vare=\tau^{(k+t)^2}3^{-k}$; we will show that $B(k,t)$ has the strong transversal property with STP-coefficient $\vare$.

Thus, let $\mathcal{B}=(B_1\LL B_{k+t})$ be an $\vare$-coherent
blockade in a
graph $G$.
We must show that 
there is a $\mathcal{B}$-transversal copy of $B(k,t)$ in $G$. Let $J$ be an optimal $\tau$-covering digraph 
for $\mathcal{B}$, with core $(A_1\LL A_{k+t})$.
Let $z=|E(J)|$.
By \ref{morecoveringdig}, every vertex of $J$ has outdegree at least one.

Since $|A_i|\ge \tau^{z} |B_i|$ for $1\le i\le k+t$,
and $\tau^{z}\ge \tau^{(k+t)^2}$,
and $\vare\le \tau^{(k+t)^2}3^{-k}$, and $\mathcal{B}$ is $\vare$-coherent, it follows 
that the blockade $(A_1\LL A_k)$ is $3^{-k}$-coherent.
\\
\\
(1) {\em We may assume that every vertex of $J$ has indegree less than $t$.}
\\
\\
Suppose that, say, $k$ is $J$-adjacent from each of $k+1\LL k+t$. Thus for each $j\in \{k+1\LL k+t\}$, there exists $X_{j,k}\subseteq B_j$
as in the definition of $\tau$-covering digraph. Since 
the blockade $(A_1\LL A_k)$ is $3^{-k}$-coherent,
by \ref{betterpath} 
there is an $(A_1\LL A_k)$-transversal $k$-vertex path in $G$ with an end-vertex in $A_k$, say with vertices $p_1\LL p_k$ in order,
where $p_i\in A_i$ for $1\le i\le k$. For each $j\in \{k+1\LL k+t\}$, choose $q_j\in X_{j,k}$ adjacent to $p_k$ (this exists, since
$X_{j,k}$ covers $A_k$). Then $q_{k+1}\LL q_{k+t}$ are pairwise nonadjacent, and nonadjacent to $p_1\LL p_{k-1}$, from the properties of
the sets $X_{j,k}$. Hence the subgraph induced on $\{p_1\LL p_k,q_1\LL q_t\}$ is a $\mathcal{B}$-transversal copy of $B(k,t)$, as required.
This proves (1).
\\
\\
(2) {\em There is a subset of $2^t+1$ elements of $\{1\LL k+t\}$, pairwise nonadjacent in $J$, and such that no two of them have
a common out-neighbour in $J$.}
\\
\\
From (1) and averaging, there is a vertex $i$ of $J$ with outdegree less than $t$; and so the set of vertices that are either equal to $i$, 
$J$-adjacent to $i$,
$J$-adjacent from $i$, or share a $J$-outneighbour with $i$, has cardinality at most $1+2(t-1)+(t-1)(t-2)=t^2-t+1$. By deleting this 
set from $J$, we obtain some digraph; again we find a vertex with outdegree at most $t-1$ in this digraph, and again
delete the corresponding set of vertices, and continue. We can repeat this at least $2^t+1$ times, since $|J|=k+t>2^t(t^2-t+1)$.
Thus we construct a set of $2^t+1$ vertices of $J$ satisfying (2). This proves (2).

\bigskip

From (2) we may assume that $\{1\LL 2^t+1\}$ is a stable set in $J$ and no two of its members have a common $J$-outneighbour.
If  $F$ is a $\mathcal{B}$-rainbow copy of $S_{t+1}$, then it has a vertex of degree $t+1$, that belongs to a block $B_i$
say, and $t+1$ other vertices, in blocks $B_j\;(j\in I)$ say, where $I\subseteq \{1\LL k+t\}\setminus \{i\}$ with $|I|=t+1$. 
Let us call
$(i,I)$ the {\em type} of $F$.
\\
\\
(3) {\em We may assume that for some $n\ge 0$ there are copies $F_1\LL F_n$ of $S_{t+1}$, each $(A_1\LL A_{t+2})$-transversal and  
pairwise disjoint, and all with type $(1,\{2\LL t+2\})$, and there exist $r\in \{1\LL t+2\}$ and $s\in \{t+3\LL t+k\}$
 with $|N(F\cap A_r)\cap A_s|\ge |A_s|/(2t+4)$, where $F=V(F_1)\cupcup V(F_n)$.}
\\
\\
Choose a maximal set $\mathcal{F}$ of pairwise disjoint $(A_1\LL A_{2^t+1})$-rainbow copies of $S_{t+1}$, and for $1\le i\le 2^t+1$
let $D_i$ be the set
of vertices of $A_i$ that are not in any member of $\mathcal{F}$. 
There is no $(D_1\LL D_{2^t+1})$-rainbow copy of $S_{t+1}$; and so by \ref{findstar}, the blockade $(D_1\LL D_{2^t+1})$ 
is not $3^{-2^t-1}$-coherent
(or one of the sets $D_i$ is empty); and it follows that $|D_r|< \vare 3^{2^t+1}|B_r|$
for some $r\in \{1\LL 2^t+1\}$. Since $|A_r|\ge \tau^{(k+t)^2} |B_r|$ and $\vare 3^{2^t+1}\le \tau^{(k+t)^2}/2$,
it follows that
$$|A_r\setminus D_r|\ge \left(\tau^{(k+t)^2} - \vare 3^{2^t+1}\right)|B_r|\ge \tau^{(k+t)^2} |B_r|/2.$$
Hence there is a subset $\mathcal{F}'\subseteq \mathcal{F}$ with cardinality at least 
$\tau^{(k+t)^2} |B_r|/2$, such that each of its members has a vertex in $A_r$.

There are only at most $(2^t+1)^{t+2}$ possible types of $\mathcal{A}$-rainbow copies of $S_{t+1}$; so within our set $\mathcal{F}'$ 
there is a subset of 
$\tau^{(k+t)^2} (2^t+1)^{-t-2}|B_r|/2$ of them all with the same type, and from the symmetry we may assume this common type is
$(1,\{2\LL t+2\})$, and $r\in \{1\LL t+2\}$.

Thus there are $n$ pairwise disjoint $(A_1\LL A_{t+2})$-rainbow copies of $S_{t+1}$, all of type $(1,\{2\LL t+2\})$, say $F_1\LL F_n$,
where $n\ge \tau^{(k+t)^2} (2^t+1)^{-t-2}|B_r|/2$. Let $F=V(F_1)\cupcup V(F_n)$. 
Since $\tau^{(k+t)^2} (2^t+1)^{-t-2}/2\ge \vare$, fewer than 
$\vare|B_{t+3}|$ vertices in $B_{t+3}$ have no neighbour in $F\cap A_r$, and so 
$$|N(F\cap A_r)\cap A_{t+3}|\ge |A_{t+3}|-\vare|B_{t+3}|\ge |A_{t+3}|/2\ge  |A_{t+3}|/(2t+4).$$
We deduce that (3) holds, setting $s=t+3$. This proves (3).

\bigskip
Let us choose $n$ minimum satisfying (3), and 
let $F=V(F_1)\cupcup V(F_n)$. We recall that $z=|E(J)|$.
\\
\\
(4) {\em For each $i\in \{1\LL t+2\}$, there is no $j\in \{1\LL k+t\}\setminus \{i\}$ such that 
$$|A_j\setminus N(F\cap A_i)|\le \tau^{z+1}|B_j|.$$ 
Consequently, if $j\in \{1\LL k+t\}\setminus \{i\}$
and 
$|N(F\cap A_i)\cap A_j|\ge \tau^{z+1}|B_j|$,
then $j$ is $J$-adjacent from $i$, and so $j\ge t+3$.}
\\
\\
For the first assertion, suppose that such $i,j$ exist. It follows that 
$$|A_j\cap N(F\cap A_i)|>|A_j|-\tau^{z+1}|B_j|.$$
From the minimality of $n$, we deduce that 
$$|A_j\cap N(F\cap A_i)|\le |A_j|/(2t+4)+\vare|B_j|.$$
Consequently
$$|A_j|-\tau^{z+1}|B_j|< |A_j|/(2t+4)+\vare|B_j|,$$
and so 
$$(1-1/(2t+4))|A_j|< (\tau^{z+1}+\vare)|B_j|.$$
Since $|A_j|\ge \tau^{z}|B_j|$, it follows that $(1-1/(2t+4))\tau^{z}<\tau^{z+1}+\vare,$
contrary to the choice of $\vare$. This proves the first assertion. The second follows immediately from \ref{coveringdig}, and so this
proves (4).

\bigskip
Choose $r,s$ as in (3). 
Let us say a vertex $v\in A_s$ is {\em good} if $v$ has a $G$-neighbour in $F\cap A_r$ and has no 
$G$-neighbour in $F\setminus A_r$. Let $C_s$ be the set of all good vertices in $A_s$.
\\
\\
(5) {\em $|C_s|\ge |A_s|/(4t+8)$.}
\\
\\
We have seen that 
$s$ is $J$-adjacent from $r$. Since no two of $1\LL t+2$ have a common $J$-outneighbour, $s$ is not $J$-adjacent from any element
in $\{1\LL t+2\}\setminus \{r\}$. From the second assertion of (4), $|N(F\cap A_i)\cap A_s|< \tau^{z+1}|B_s|$ 
for each $i\in \{1\LL t+2\}\setminus \{r\}$. Hence at most $(t+1)\tau^{z+1}|B_s|$ vertices in $B_s$ have a $G$-neighbour in 
$F\setminus A_r$;
and so at least 
$|A_s|/(2t+4)-(t+1)\tau^{z+1}|B_s|$ vertices in $A_s$ are good.
Since 
$$(t+1)\tau^{z+1}|B_s|\le \tau^{z}|B_s|/(4t+8)\le |A_s|/(4t+8)$$
we deduce that $|C_s|\ge |A_s|/(4t+8)$. This proves (5).

\bigskip
We know that $r\in \{1\LL t+2\}$, but the argument to come depends on whether $r=1$ or not. If $r=1$ let $t'=t+1$ and $F'=F\setminus A_{t+2}$, and otherwise let $t'=
t+2$ and $F'=F$. For $1\le m\le n$ let $F_m'=F_m\setminus A_{t+2}$ if $r=1$, and otherwise let $F_m'=F_m$. Thus if $r=1$, $F_m'$
is a copy of $S_t$, and otherwise $F_m'$ is a copy of $S_{t+1}$.
For each $j\in \{t'+1\LL k+t\}$ with $j\ne s$, let $C_{j}$ be the set of all vertices in $A_{j}$ with no neighbour in $F'$. 
\\
\\
(6) {\em $|C_j|\ge  |A_{j}|/(4t+8)$ for each $j\in \{t'+1\LL k+t\}$.}
\\
\\
Let $j\in \{t'+1\LL k+t\}$. By (5) we may assume that $j\ne s$.
From the 
minimality of $n$, for $1\le i\le t'$ at most $|A_{j}|/(2t+4)+\vare|B_{j}|$ vertices in $A_{j}$ have a neighbour in $F'\cap A_i$; and so
at most $(t+2)(|A_{j}|/(2t+4)+\vare|B_{j}|)$ have a neighbour in $F$. Thus 
$$|C_{j}|\ge |A_{j}|-|A_{j}|/2-(t+2)\vare|B_{j}|= |A_{j}|/2- (t+2)\vare|B_{j}|\ge |A_{j}|/(4t+8).$$
This proves (6).

\bigskip
From (6) the blockade $(C_{t'+1}\LL C_{k+t})$ is $(4t+8)\vare\tau^{-z}$-coherent, 
and since
$(4t+8)\vare\tau^{-z}\le 3^{1-(k+t)+t'}$ (because $z\le (k+t)(k+t-1)$ and $t'-t\le 2$)
\ref{betterpath} implies that there is a $(C_{t'+1}\LL C_{k+t})$-transversal path $P$ with end-vertex in $C_s$.
Let its end in $C_s$ be $v$. From the definition of $C_s$, there exists $u\in F'\cap A_r$ adjacent to $v$. Let $u\in F_m'$ say, where
$1\le m\le n$. Since $v\in C_s$, $v$ is good and so has no neighbour in $F_m'$ except $u$; and since every other vertex of $P$
belongs to some set $C_j$ where $j\in \{t'+1\LL k+t\}\setminus \{s\}$, the edge $uv$ is the only edge of $G$ between $V(P)$
and $V(F_m')$. Thus the union of $P,F_m$ and the edge $uv$ is a $\mathcal{B}$-rainbow copy of $B(k,t)$. This proves \ref{brooms}.~\bbox

Thus, if $t\le \log_2 k-(2+o(1))\log\log_2 k $, \ref{brooms} tells us that $B(k,t)$ has the strong transversal property, and for 
$t>(1+o(1))\log_2 k $, \ref{betterclaw} tells us that it does not. More exactly, for $t\ge 2$, if $k \ge 2^t(t^2-t+1)-t+1$ then $B(k,t)$
has the property, and if $2\le k\le 2^t-t$ then it does not. We have not decided the values of $k$ in the  middle, except when $t=2$.
In that case \ref{brooms} tells us that $B(k,2)$ has the property when $k\ge 11$, but this can be improved to:
\begin{thm}\label{2broom}
If $k=1$ or $k\ge 3$ then $B(k,2)$ has the strong transversal property, and if $k=2$ it does not.
\end{thm}
\Proof
We just sketch the proof, since it is similar to that of \ref{brooms}. The claims for $k=1$ and $k=2$ follow from \ref{path}
and \ref{noclaw} respectively, so we assume that $k\ge 3$. With an appropriate choice of $\tau$ and $\vare$, we choose an optimal
$\tau$-covering digraph $J$; and we may assume no vertex has indegree more than one in $J$, as in the proof of \ref{brooms}, and every
vertex has outdegree at least one, by \ref{morecoveringdig}. Consequently $J$ is a disjoint union of directed cycles. 
Next we use a lemma (we omit the proof here), that for any $\tau$-covering digraph $J$ and for every directed cycle of $J$, some edge of the cycle
is in a directed cycle of length two (this is true in general, not just in the present context). Consequently $J$ is the disjoint
union of directed cycles of length two, and in particular, $|J|$ is even. Thus there are three pairwise nonadjacent
vertices of $J$, say $1,2,3$; and as in the proof of \ref{brooms} we find many
pairwise disjoint copies of $S_2$, all $(A_1,A_2,A_3)$-rainbow and all with middle vertex in $A_1$. (Note that we look for copies of
$S_2$, not for copies of $S_3$, which is what setting $t=2$ in the proof of \ref{brooms} would suggest.) We may 
assume that $4,5$ and $6$ are $J$-adjacent to and from $1,2,3$ respectively. Now the proof is 
finished more-or-less as in \ref{brooms}; with notation as in \ref{brooms}, if $r=1$ we follow the proof of \ref{brooms}, that is,
we delete from $A_4\LL A_{k+2}$ the small number of vertices with the wrong adjacency to $F$, and then apply \ref{betterpath}
to the resulting blockade $(C_4,C_5\LL C_{k+2})$, finding a path with first vertex in $C_4$. If $r=2$ say, we do the same,
but apply \ref{betterpath}
to the blockade $(C_5\LL C_{k+2})$, finding a path with first vertex in $C_5$, and then turning this into a copy of
$B(2,k)$ by adding $F_m$ and a vertex of $X_{4,1}$. We omit further details.~\bbox

\section{The cycle}

Our remaining results all concern looking for an anticomplete pair of sets that have polynomial
size rather than linear size. So, we are not working with the strong transversal property any more, nor with $\vare$-coherence.

In this section we prove \ref{cycle}. 
We handle the cases $k=4$ and $k\ge 5$ separately.
Both proofs are related to the proofs of theorems
in~\cite{polyC5}. (These theorems have recently been superceded by theorems in~\cite{fivehole}, but the proof methods 
of the 
latter are quite different.)

We will need the following lemma:
\begin{thm}\label{manyedges}
Let $0<\vare\le 1/2$, and let $\mathcal{B}=(B_1,B_2)$ be a
blockade in a
graph $G$, with local degree less than $\vare W$ and $(\vare W,\vare W^{c})$-cohesive where $W$ is its width.
Then $\vare W^{c}> 1$; and if $X\subseteq B_1$ with $|X|\ge 2\vare W$, there are fewer than $\vare W^{c}$ vertices in $B_2$
that have 
at most $W^{1-c}/2$ neighbours in $X$. 
\end{thm}
\Proof
Let $v\in B_1$. Since $\mathcal{B}$ has local degree less than $\vare W$, $v$ has at most $\vare W$ neighbours in $B_2$; and so has
at least $(1-\vare)W\ge \vare W$ non-neighbours in $B_2$. Thus $\mathcal{B}$ is not $(\vare W,1)$-cohesive, and so $\vare W^{c}>1$.
This proves the first assertion.

Suppose the second assertion is false; then there exists $Y\subseteq B_2$ with $|Y|=\lceil \vare W^{c}\rceil$, such that
every vertex in $Y$ has at most $W^{1-c}/2$ neighbours in $B_2$. Since $\vare W^{c}> 1$, it follows that
$|Y|\le 2\vare W^{c}$.
Hence at most $|Y| W^{1-c}/2\le \vare W$ vertices in $X$ have a neighbour in $Y$, and since $|X|\ge 2\vare W$, $X$ has a subset
of cardinality at least $\vare W$ that is anticomplete to $Y$, a contradiction. This proves \ref{manyedges}.~\bbox

First we show the following, which is a strengthening of the $k=4$ case of \ref{cycle}:

\begin{thm}\label{4cycle}
Let $\vare=1/4$, and let $\mathcal{B}=(B_1\LL B_4)$ be a
blockade in a
graph $G$, with local degree less than $\vare W$ and $(\vare W,\vare W^{1/3})$-cohesive where $W$ is its width.
Then there is a $\mathcal{B}$-transversal copy of a cycle of length four in $G$.
\end{thm}
\Proof 
Let $\mathcal{B}=(B_1\LL B_4)$ be a
blockade in a
graph $G$, with local degree less than $\vare W$ and $(\vare W,\vare W^{1/3})$-cohesive where $W$ is its width.
From \ref{manyedges}:
\\
\\
(1) {\em $\vare W^{1/3}> 1$; and if $i,j\in \{1\LL 4\}$ are distinct, and $X\subseteq B_i$ with $|X|\ge 2\vare W$, and
$Y$ is a set of vertices in $B_j$ each with at most $W^{2/3}/2$ neighbours in $X$, then $|Y|< \vare W^{1/3}$.}

\bigskip
Let $v_3\in B_3$ and $v_4\in B_4$ be adjacent. We say the edge $v_3v_4$ is
\begin{itemize}
\item {\em 1-good} if $v_3$ has at least $W^{2/3}/2$ neighbours in $B_4$;
\item {\em 2-good} if it is 1-good and $v_4$ has at least $W^{2/3}/2$ neighbours in $B_2$ that are nonadjacent to $v_3$.
\end{itemize}
We claim:
\\
\\
(2) {\em More than half the edges between $B_3,B_4$ are 2-good.}
\\
\\
By (1), fewer than $\vare W^{1/3}$ vertices in $B_3$ have at most $W^{2/3}/2$ neighbours in $B_4$, so at most
$\vare W^{1/3}W^{2/3}/2=\vare W/2$ edges between $B_3,B_4$ are not 1-good. Now let $v_3\in B_3$ and let $N_4$
be the set of its neighbours in $B_4$, and let $N_2$ be the set of its neighbours in $B_2$. Thus $|N_2|<\vare |B_2|$, and so 
$|B_2\setminus N_2|\ge (1-\vare) |B_2|\ge 2\vare W$; so by (1), fewer than $\vare W^{1/3}$ vertices in $N_4$
have fewer than $W^{2/3}/2$ neighbours in $B_2\setminus N_2$. Consequently at most $\vare W^{1/3}$ of the edges between
$v_3$ and $N_4$ are 1-good and not 2-good. Since this holds for every choice of $v_3\in B_3$, it follows that
at most $\vare W^{1/3}|B_3| $ edges between $B_3,B_4$ are 1-good and not 2-good. Hence in total, at most
$\vare W^{1/3}|B_3| +\vare W/2$  edges between $B_3,B_4$ are not 2-good. But at least $|B_3|-\vare W^{1/3}$ vertices in
$B_3$ have at least $W^{2/3}/2$ neighbours in $B_4$, so there are at least 
$$(|B_3|-\vare W^{1/3})W^{2/3}/2$$
edges between $B_3,B_4$. Since
$$\vare W^{1/3}|B_3| +\vare W/2< (|B_3|-\vare W^{1/3})W^{2/3}/4$$
because, for instance, 
$$(W^{2/3}/4-\vare W^{1/3})|B_3|\ge (W^{2/3}/4-\vare W^{1/3})W> \vare W/2+\vare W/4,$$
this proves (2).
\\
\\
(3) {\em We may assume that there exist anticomplete subsets $C_1\subseteq B_1$ and $C_2\subseteq B_2$, with 
$|C_1|,|C_2|\ge W^{2/3}/2$.}
\\
\\
By (2), and the same statement with $B_3,B_4$ exchanged and $B_1,B_2$ exchanged, it follows that there is an edge $v_3v_4$
with $v_3\in B_3$ and $v_4\in B_4$, such that $v_4$ has at least $W^{2/3}/2$ neighbours in $B_2$ that are nonadjacent to $v_3$,
and $v_3$ has at least $W^{2/3}/2$ neighbours in $B_1$ that are nonadjacent to $v_4$. Let $C_1$ be the set of vertices in
$B_1$ that are adjacent to $v_3$ and not to $v_4$, and define $C_2\subseteq B_2$ similarly. If there is an edge
between $C_1,C_2$, then adding $v_3,v_4$ makes a $\mathcal{B}$-transversal cycle of length four; so we may assume there
is no such edge. This proves (3).

\bigskip
Choose $C_1,C_2$ as in (3). 
\\
\\
(4) {\em There exists $v_3\in B_3$ with at least $\vare W^{1/3}$ neighbours in $C_1$ and at least $\vare W^{1/3}$ 
neighbours in $C_2$.}
\\
\\
Suppose that there is a set $A_3\subseteq B_3$ with cardinality $\lceil 2\vare W\rceil$, such that each of its members
has fewer than $\vare W^{1/3}$ neighbours 
in $C_1$.
By (1), fewer than  $ \vare W^{1/3}$ vertices in $C_1$ have fewer than $W^{2/3}/2$ neighbours in $A_3$;
and so there are at least 
$$(|C_1|-\vare W^{1/3})W^{2/3}/2$$
edges between $C_1, A_3$. But from the definition of $A_3$, there are at most
$|A_3|\vare W^{1/3}$ such edges; so
$$(|C_1|-\vare W^{1/3})(W^{2/3}/2) <|A_3|(\vare W^{1/3}).$$
Since $|C_1|\ge W^{2/3}/2$, and 
$|A_3|\le 3\vare W$ (since $\vare=1/4$ and $W\ge 64$), it follows that
$$(W^{2/3}/2-\vare W^{1/3})(W^{2/3}/2) < 3\vare W(\vare W^{1/3}),$$
which (since $\vare=1/4$) simplifies to 
$W^{1/3} < 2$, a contradiction since $W^{1/3}\ge \vare^{-1}=4$.
Thus there are fewer than $2\vare W$ vertices in $B_3$ with 
fewer than $\vare W^{1/3}$ neighbours            
in $C_1$; and fewer than $2\vare W$ vertices in $B_3$ with fewer than $\vare W^{1/3}$ neighbours                
in $C_2$, similarly. Since $4\vare W=W\le |B_3|$, 
this proves (4).

\bigskip
Choose $v_3$ as in (4), and for $i = 1,2$ let $A_i$ be the set of neighbours of $v_3$ in $C_1$ and in $C_2$
respectively. Since $\mathcal{B}$ is $(\vare W,\vare W^{1/3})$-cohesive, fewer than $\vare W$ vertices in $B_4$
have no neighbour in $A_1$; fewer than $\vare W$ vertices in $B_4$
have no neighbour in $A_2$; and fewer than $\vare W$ vertices in $B_4$
are adjacent to $v_3$. Since $3\vare W<W\le |B_4|$, there exists $v_4\in B_4$ with a neighbour $v_1\in A_1$, and a neighbour 
$v_2\in A_2$,
and non-adjacent to $v_3$. But then there is a $\mathcal{B}$-transversal 4-cycle induced on $\{v_1,v_2,v_3,v_4\}$.
This proves \ref{4cycle}.~\bbox

To complete the proof of \ref{cycle}, next we prove the following, a strengthening of \ref{cycle} when $k\ge 5$:
\begin{thm}\label{kcycle}
Let $k\ge 5$ be an integer, and let $\vare=1/(3k)$.
Let $\mathcal{B}=(B_1\LL B_k)$ be a
blockade in a
graph $G$, with local degree less than $\vare W$ and $(\vare W,\vare W^{1/2})$-cohesive where $W$ is its width.
Then there is a $\mathcal{B}$-transversal copy of a cycle of length $k$ in $G$.
\end{thm}
\Proof
Let $\mathcal{B}=(B_1\LL B_k)$ be a
blockade in a
graph $G$, with local degree less than $\vare W$ and $(\vare W,\vare W^{1/2})$-cohesive where $W$ is its width.
From \ref{manyedges}:
\\
\\
(1) {\em $\vare W^{1/2}>1$; and if $i,j\in \{1\LL k\}$ are distinct, and $X\subseteq B_i$ with $|X|\ge 2\vare W$, and
$Y$ is a set of vertices in $B_j$ each with fewer than $W^{1/2}/2$ neighbours in $X$, then $|Y|< \vare W^{1/2}$.}

\bigskip

Next we prove the following:
\\
\\
(2) {\em Let $i_1,i_2\LL i_s\in \{1\LL k\}$ be distinct, and let $D_{i_r}\subseteq B_{i_r}$ for $1\le r\le s$, such that
$|D_{i_1}|\ge \vare W^{1/2}$, and $|D_{i_r}|\ge \vare s |B_{i_r}|$ for $2\le r\le s$. There there is an induced path of $G$
with vertices $v_{i_1}\CC v_{i_s}$, where $v_{i_r}\in D_{i_r}$ for $1\le r\le s$.}
\\
\\
We proceed by induction on $s$; if $s=1$ the result is trivial, so we assume that $s\ge 2$ and the result holds for $s-1$.
Let  $i_1,i_2\LL i_s$ and $D_{i_r}\subseteq B_{i_r}$ for $1\le r\le s$ as above.
From the symmetry we may assume that $i_r=r$ for $1\le r\le s$.
Since $s\ge 2$, and consequently $|D_{2}|\ge s\vare |B_{2}|\ge 2\vare|B_{2}|$,  (1) implies
that there exists $v_{1}\in D_{1}$ with at least $W^{1/2}/2$ neighbours in $D_{2}$. Let $E_{2}$ be the set of these 
neighbours, and for $3\le r\le s$ let $E_{r}$ be the set of vertices in $D_{r}$ nonadjacent to $v_{1}$. Since $v_{1}$
has fewer than $\vare |B_{r}|$ neighbours in $B_{r}$, it follows that 
$|E_{r}|\ge |D_{r}|-\vare |B_{r}|\ge \vare(s-1)|B_{r}|$. Hence from the inductive hypothesis applied to $2\LL s$
and the sets $E_{2}\LL E_{s}$, there is an induced path of $G$
with vertices $v_{2}\CC v_{s}$, where $v_{r}\in E_{r}$ for $2\le r\le s$. Adding $v_{1}$ and the edge
$v_{1}v_2$ gives a path satisfying (2). This proves (2).

\bigskip
From (1), all vertices in $B_1$ except at most $\vare W^{1/2}$ have at least $W^{1/2}/2$ neighbours in $B_2$, and the same 
for $B_3$; so there exists $v_1\in B_1$ with at least $W^{1/2}/2$ neighbours in $B_2$ and at least $W^{1/2}/2$ 
neighbours in $B_3$. For $i = 2,3$ let $A_i$ be the set of neighbours of $v_1$ in $B_i$, and for $4\le i\le k$ let $A_i$
be the set of vertices in $B_i$ that are nonadjacent to $v_1$. Thus $|A_2|,|A_3|\ge W^{1/2}/2$. Since
$v_1$ has at most $\vare|B_j|$ neighbours in $B_j$, it follows that $|A_j|\ge (1-\vare)|B_j|$ for $4\le j\le k$.

All except at most $\vare W$ vertices in $A_4$ have a neighbour in $A_2$, so we may choose
$C_2\subseteq A_2$ minimal such that for some $j\in \{4\LL k\}$, at least 
$|B_j|/3$ vertices in $A_j$ have a neighbour in $C_2$. Choose some such $j$; and from the symmetry we may assume that $j=4$.
Let $C_4$ be the set of vertices in $A_4$
that have a neighbour in $C_2$. For $5\le i\le k$, let $C_i$ be the set of vertices in $A_i$ with no 
neighbour in $C_2$. Thus $|C_4|\ge |B_4|/3$; and from the minimality of $C_2$, it follows that fewer than $|B_i|/3+\vare|B_i|$
have a neighbour in $C_2$, for each $i\in \{5\LL k\}$.
Hence for each $i\in \{5\LL k\}$, 
$$|C_i|\ge |A_i|-|B_i|/3-\vare|B_i|\ge (1-\vare)|B_i|-|B_i|/3-\vare|B_i|=(2/3-2\vare)|B_i|\ge |B_i|/3.$$
By (1), at most $\vare W^{1/2}$
vertices in $A_3$ have fewer than $W^{1/2}/2$ neighbours in $C_5$ (this is where we use $k\ge 5$); 
and since $|A_3|\ge W^{1/2}/2 >\vare W^{1/2}$, there exists
$v_3\in A_3$ with at least $W^{1/2}/2$ neighbours in $C_5$.

Let $D_5$ be the set of neighbours of $v_3$ in $C_5$; let $C_1=B_1$, and for
$i\in \{1,4\}\cup \{6\LL k\}$, let $D_h$ be the set of vertices in $C_h$ nonadjacent to $v_3$. Hence
$|D_5|\ge W^{1/2}/2$, and 
$$|D_h|\ge |C_h|-\vare|B_h|\ge (2/3-3\vare)|B_i|$$
for $i\in \{1,4\}\cup \{6\LL k\}$.

Every vertex in $D_4$ has a neighbour in $C_2$, which may or may not be adjacent to $v_3$. So either at least $|D_4|/2$
vertices in $D_4$ have a neighbour in $C_2$ nonadjacent to $v_3$, or at least $|D_4|/2$
vertices in $D_4$ have a neighbour in $C_2$ adjacent to $v_3$. We handle these two cases separately.

First, assume that at least $|D_4|/2$
vertices in $D_4$ have a neighbour in $C_2$ nonadjacent to $v_3$; let $D_2$ be the set of vertices in $C_2$ nonadjacent to $v_3$,
and let $D_4'$ be the set of vertices in $D_4$ with a neighbour in $D_2$. Thus 
$$|D_4'|\ge |D_4|/2\ge (2/3-3\vare)|B_4|/2\ge (k-3)\vare |B_4|,$$
since $\vare\le 1/(3k)$.
By (2), 
there is an induced path $P$ of $G$
$v_{5}\DD v_{6}\CC v_k\DD v_4$, where $v_{r}\in D_{r}$ for $5\le r\le k$ and $v_{4}\in D_4'$.
Choose $v_2\in D_2$ adjacent to $v_4$; then the union of $P$ and the path 
$v_4\DD v_2\DD v_1\DD v_3\DD v_5$ gives a $\mathcal{B}$-transversal cycle satisfying the theorem.

Now we assume that  at least $|D_4|/2$
vertices in $D_4$ have a neighbour in $C_2$ adjacent to $v_3$; let $D_2$ be a subset of $C_2$, all
adjacent to $v_3$, minimal such that either at least $|D_4|/2$ vertices in $D_4$, or at least $|D_1|/2$ vertices in $D_1$,
have a neighbour in $D_2$. 

Suppose there is a set $D_4'\subseteq D_4$ with cardinality at least $|D_4|/2$, all with a neighbour in $D_2$. From the minimality of $D_2$, at most
$|D_1|/2 +\vare|B_1|$ vertices in $D_1$ have a neighbour in $D_2$, and so there is a subset $D_1'\subseteq D_1$
with cardinality at least $|D_1|/2-\vare|B_1|\ge (k-2)\vare |B_1|$,  anticomplete to $D_2$. 
By (2),
there is an induced path $P$ of $G$
with vertices $v_5\CC v_k\DD v_1'\DD v_4$, where $v_r\in D_{r}$ for $5\le r\le k$ and $v_{1}'\in D_1'$ and $v_4\in D_4'$.
Choose $v_2\in D_2$ adjacent to $v_4$; then the union of $P$ and the path with vertices
$v_4\DD v_2\DD v_3\DD v_5$ gives a $\mathcal{B}$-transversal cycle satisfying the theorem.

Finally we may assume that there is a set $D_1'\subseteq D_1$ with cardinality at least $|D_1|/2$, all with a neighbour in $D_2$.
From the minimality of $D_2$, at most
$|D_4|/2 +\vare|B_4|$ vertices in $D_4$ have a neighbour in $D_2$, and so there is a subset $D_4'\subseteq D_4$
with cardinality at least $|D_4|/2-\vare|B_4|\ge (k-2)\vare |B_4|$,  anticomplete to $D_2$. 
By (2),
there is an induced path $P$ of $G$
with vertices $v_5\CC v_k\DD v_4\DD v_1'$, where $v_r\in D_{r}$ for $5\le r\le k$ and $v_4\in D_4'$ and $v_{1}'\in D_1'$.
Choose $v_2\in D_2$ adjacent to $v_1'$; then the union of $P$ and the path with vertices
$v_1'\DD v_2\DD v_3\DD v_5$ gives a $\mathcal{B}$-transversal cycle satisfying the theorem. This proves \ref{kcycle}.~\bbox

We do not know whether the exponents of $1/3$ (in \ref{4cycle}) and $1/2$ (in \ref{kcycle}) are best possible.

\section{Ordered transversal subgraphs}

Now we turn to excluding ordered graphs.
We begin with \ref{tree}, which we restate:

\begin{thm}\label{tree2}
If $H$ is an ordered tree with $k\ge 1$ vertices, then there exists $\vare>0$ with the following property.
Let $\mathcal{B}=(B_1\LL B_k)$ be a
blockade in a
graph $G$, with local degree less than $\vare W$ and $(\vare W,\vare W^{1/(k-1)})$-cohesive where $W$ is its width.
Then there is an ordered $\mathcal{B}$-transversal copy of $H$.
\end{thm}

To prove this we need to prove a strengthening (which implies \ref{tree2} by setting $c=1/(k-1)$):
\begin{thm}\label{counttree}
Let $H$ be an ordered tree with $k\ge 1$ vertices, and let $c>0$ with $(k-1)c\le 1$. Let $\vare=4^{1-k}$.
Let $\mathcal{B}=(B_1\LL B_k)$ be a
blockade in a
graph $G$, with local degree less than $\vare W$ and $(\vare W,\vare W^{c})$-cohesive where $W$ is its width.
Then there are at least $4^{1-k}W^{k-(k-1)c}$ ordered $\mathcal{B}$-transversal copies of $H$.
\end{thm}
\Proof
We proceed by induction on $k$. The result is trivial for $k=1$, so we assume that $k\ge 2$ and the result holds for $k-1$.
Let the ordering of $H$ be $p_1\LL p_k$. We may assume that $p_k$ has degree one in $H$, and $p_{k-1}$ is its unique neighbour.
Let $\mathcal{B}=(B_1\LL B_k)$ be a
blockade in a
graph $G$, with local degree less than $\vare W$ and $(\vare W,\vare W^{c})$-cohesive where $W$ is its width.
We may assume that $|B_i|=W$ for $1\le i\le k$.
From \ref{manyedges}:
\\
\\
(1) {\em $\vare W^{c}> 1$; and if $X\subseteq B_{k}$ with $|X|\ge 2\vare W$, there are fewer than $\vare W^{c}$ vertices in $B_{k-1}$
that have
at most $W^{1-c}/2$ neighbours in $X$.}

\bigskip
In particular, there are at least $(W-\vare W^{c})W^{1-c}/2\ge W^{2-c}/4$ edges between $B_{k-1}$ and $B_k$, so if $k=2$ the result is true. Thus we may assume that $k\ge 3$.

Let $H'$
be the ordered tree obtained from $H$ by deleting $p_k$ (with ordering $p_1\LL p_{k-1}$), and similarly let $H''$
be the ordered forest obtained by deleting both $p_k, p_{k-1}$. Let $\mathcal{B}'=(B_1\LL B_{k-1})$, and $\mathcal{B}''=(B_1\LL B_{k-2})$.
Let $\mathcal{H}$ be the set of all  ordered $\mathcal{B}$-transversal copies of $H$,
let $\mathcal{H}'$ be the set of all ordered $\mathcal{B}'$-transversal copies of $H'$, and let $\mathcal{H}''$ be the set of all ordered $\mathcal{B}''$-transversal copies of $H''$.

For each $F\in \mathcal{H}''$,  let $n(F)$ be 
the number of vertices $v\in B_{k-1}$ such that the
subgraph induced on $V(F)\cup \{v\}$ is an ordered $\mathcal{B}'$-transversal copy of $H'$. 
Let $m(F)$ be the number of edges $uv$ with $u\in B_k$ and $v\in B_{k-1}$ such that 
the
subgraph induced on $V(F)\cup \{u,v\}$ is an ordered $\mathcal{B}$-transversal copy of $H$. 
\\
\\
(2) {\em For each $F\in \mathcal{H}''$, $m(F)\ge (n(F)-\vare W^c)W^{1-c}/2$.}
\\
\\
Let $F\in \mathcal{H}''$ and let $N$ be the set of vertices $v\in B_{k-1}$ such that the
subgraph induced on $V(F)\cup \{v\}$ is an ordered $\mathcal{B}'$-transversal copy of $H'$. Let $X$ be the set of vertices
in $B_k$ with no neighbours in $V(F)$. Thus 
$$|X|\ge |B_1|-(k-2)\vare|B_1|\ge 2\vare W,$$
(since $\vare\le 1/k$),
and so by (1), there are fewer than  $\vare W^{c}$ vertices in $N$
that have
at most $W^{1-c}/2$ neighbours in $X$. All the others have more than $W^{1-c}/2$ neighbours in $X$, and every such edge
contributes to $m(F)$. This proves (2).

\bigskip
Summing $n(F), m(F)$ and $1$ over all $F\in \mathcal{H}''$ gives $|\mathcal{H}|, |\mathcal{H}'|$ and $|\mathcal{H}''|$ 
respectively, so 
by summing the inequality of (2) over all $F\in \mathcal{H}''$, we deduce that
$$2|\mathcal{H}|\ge W^{1-c} |\mathcal{H}'|  - \vare W|\mathcal{H}''|\ge W^{1-c} |\mathcal{H}'|  - \vare W^{k-1}.$$
But from the inductive hypothesis, 
$$|\mathcal{H}'|\ge 4^{2-k}W^{k-1-(k-2)c}\ge 2\vare W^{k-2+c}$$
(the latter since $c\le 1/(k-1)$ and $\vare= 4^{1-k}$).
Consequently
$$2|\mathcal{H}|\ge W^{1-c} |\mathcal{H}'|/2\ge 4^{2-k}W^{k-1-(k-2)c} W^{1-c}/2=2\cdot 4^{1-k}W^{k-(k-1)c}.$$
This proves \ref{counttree}.~\bbox

The exponent of $1/(k-1)$ in \ref{tree2} is best possible in the sense that for the tree $S_{k-1}$, the exponent
cannot be replaced by any larger constant, as we shall see. But perhaps it can be replaced by $1/d$ where $d$ is the 
maximum degree of the tree? We propose:

\begin{thm}\label{treeconj}
{\bf Conjecture: }If $H$ is an ordered tree with $k\ge 1$ vertices and maximum degree $d$, then there exists $\vare>0$ with the following property.
Let $\mathcal{B}=(B_1\LL B_k)$ be a
blockade in a
graph $G$, with local degree less than $\vare W$ and $(\vare W,\vare W^{1/d})$-cohesive where $W$ is its width.
Then there is an ordered $\mathcal{B}$-transversal copy of $H$.
\end{thm}
The next result shows that this is true for caterpillars:
\begin{thm}\label{orderedcaterpillar}
Let $H$ be an ordered caterpillar with $k\ge 1$ vertices and maximum degree $d$, and let $\vare=4^{-d}/k$.
Let $\mathcal{B}=(B_1\LL B_k)$ be a
blockade in a
graph $G$, with local degree less than $\vare W$ and $(\vare W,\vare W^{1/d})$-cohesive where $W$ is its width.
Then there is an ordered $\mathcal{B}$-transversal copy of $H$.
\end{thm}
For inductive purposes it is helpful to prove something stronger. If $H$ is a caterpillar, there is a path of $H$ containing
all vertices of $H$ with degree more than one. If there is such a path with one end $v$ we call $v$ a {\em head} of the caterpillar. (Thus, the head is not necessarily unique.)
We will show:
\begin{thm}\label{orderedcaterpillar2}
Let $H$ be an ordered caterpillar with $k\ge 1$ vertices, with ordering $v_1\LL v_k$ where $v_1$ is a head.
Let $v_1$ have degree $d_1$, and let every vertex of $H$ have degree at most $d$, and let $\vare=4^{-d}/k$.
Let $\mathcal{B}=(B_1\LL B_k)$ be a
blockade in a
graph $G$, with local degree less than $\vare W$ and $(\vare W,\vare W^{1/d})$-cohesive where $W$ is its width.
Let $C_1\subseteq B_1$, where $|C_1|\ge 4^{d_1-d}W^{d_1/d}$. Then there is an ordered $(C_1,B_2\LL B_k)$-transversal copy of $H$.
\end{thm}
\Proof
Let 
$G$, $\mathcal{B}=(B_1\LL B_k)$ and $C_1$ be as in the theorem. We may assume that $|B_1|\LL |B_k|=W$.
If $k=1$ the result is trivial; and if $k=2$, then $d_1=d=1$, and so $|C_1|\ge W$, and therefore there is an edge between $C_1, B_2$
and the claim holds. So we may assume that $k\ge 3$, and proceed by induction on $k$. 

Suppose first that $d_1\ge 2$, and let $v_k$ say be a neighbour of $v_1$ that has degree one in $H$.
Let $H'$ be obtained from $H$ by deleting $v_k$.
By \ref{manyedges}, there are fewer than $\vare W^{1/d}$ vertices in $C_1$
that have
at most $W^{1-1/d}/2$ neighbours in $B_k$. Hence there are at least 
$(|C_1|-\vare W^{1/d})W^{1-1/d}/2\ge |C_1| W^{1-1/d}/4$ edges between $C_1$ and $B_k$, since $|C_1|\ge W^{d_1/d}\ge 2\vare W^{1/d}$.
Consequently some vertex $u_k\in B_k$ has at least $|C_1| W^{-1/d}/4$ neighbours in $C_1$. Let $C_1'$ be the set of these neighbours;
then 
$$|C_1'|\ge |C_1| W^{-1/d}/4\ge 4^{d_1-d}W^{d_1/d}W^{-1/d}/4\ge 4^{d_1-1-d}W^{(d_1-1)/d}.$$
For $2\le i\le k-1$, let $B_i'$ be the set of vertices in $B_i$ nonadjacent to $u_k$; so $|B_i'|\ge (1-\vare)W$.
Let the blockade $(B_1,B_2'\LL B_{k-1}')$ have width $W'$ say; then $W'\ge (1-\vare)W$. Let
$\vare'=4^{-d}/(k-1)$; then $\vare' W'\ge \vare W$, and so $(B_1,B_2'\LL B_{k-1}')$
has 
local degree less than $\vare' W'$ and is $(\vare' W',\vare' (W')^{1/d})$-cohesive.
(Note that $\vare' (W')^{1/d}\ge \vare W^{1/d}$.)
From the inductive hypothesis, there is an ordered $(C_1', B_2'\LL B_{k-1}')$-transversal copy of $H\setminus \{v_k\}$;
and adding $u_k$ gives an ordered $\mathcal{B}$-transversal copy of $H$ containing
a vertex of $C_1$.

So we may assume that $d_1=1$; let $v_2$ be the unique neighbour of $v_1$ in $H$, and let $H'$ be obtained from
$H$ by deleting $v_1$. Thus $v_2$ is a head of $H'$.

By \ref{manyedges}, since $|C_1|\ge 4^{1-d}W^{1/d}>\vare W^{1/d}$, there is a vertex $u_1\in C_1$
with at least $W^{1-1/d}/2$ neighbours in $B_2$; let $C_2$ be the set of these neighbours.
For $3\le i\le k$, let $B_i'$ be the set of vertices in $B_i$ nonadjacent to $u_1$, so $|B_i'|\ge (1-\vare)W$.
Let the blockade $(B_2,B_3'\LL B_{k}')$ have width $W'$ say. Then $W'\ge (1-\vare)W$, and as before, it has
local degree less than $\vare' W'$ and is $(\vare' W',\vare' (W')^{1/d})$-cohesive where $\vare'=4^{-d}/(k-1)$.
From the inductive hypothesis, there is an ordered $(C_2, B_3'\LL B_{k}')$-transversal copy of $H\setminus \{v_1\}$,
and adding $u_1$ gives a an ordered $\mathcal{B}$-transversal copy of $H$ containing
a vertex of $C_1$. This proves \ref{orderedcaterpillar2}.~\bbox

Finally, let us see that the exponents in \ref{orderedcaterpillar} and \ref{orderedcaterpillar2}
cannot be replaced by any larger constant. 
We need the following three lemmas:
\begin{thm}\label{bigrandom1}
Let $t\ge 3$ be an integer. Let $0<\vare<1$ be rational, and let $c>d>1/t$, where $d-1/t < (c-1/t)/(t-1)$, and $d<2/t$, and $c,d$ are rational. 
Let $K>\vare  \log (e/\vare)/ (-\log (1-\vare))$.
Let $n$ be an integer such that
$n^c, n^d, n^{1-d}, n^{1/t}, \vare n$ are all integers. If $n$ is sufficiently large,
there is a graph with bipartition $A,B$ , where $|A|=n$ and $|B|=n^{2/t}$, such that
\begin{itemize}
\item every vertex in $A$ has degree at most $t-1$;
\item for every $X\subseteq A$ with $|X|\ge \vare n^c$, there are at least $Kn^d$ vertices in $B$ that have a neighbour in $X$;
\item for every $X\subseteq A$ with $|X|\ge \vare n$, there are at least $(\vare/e) |B|$ vertices in $B$ with a neighbour in $X$; and
\item every vertex in $B$ has less than $n^{1-d}$ neighbours in $A$.
\end{itemize}
\end{thm}
\Proof
Let $A,B$ be disjoint sets of cardinalities $n, n^{2/t}$ respectively. For each $v\in A$, choose $v_1\LL v_{t-1}$ in $B$ uniformly
and independently at random (and therefore not necessarily distinct), and add edges to make $v$ adjacent to $v_1\LL v_{t-1}$. Let $G$
be the graph this constructs. We claim that if $n$ is sufficiently large then with high probability $G$ satisfies the theorem.

Let $X\subseteq A$ with $|X|= \vare n^c$, and let $Y\subseteq B$ with $|Y|=Kn^d$. The probability that for every vertex in $X$, all its neighbours are in $Y$, is 
$$\left(\left(Kn^{d}\right)/|B|\right)^{\vare(t-1)n^c}=\left(Kn^{d-2/t}\right)^{\vare(t-1)n^c}.$$
By \ref{binom} there are at most $((e/\vare)n^{1-c})^{n^c}$ choices of $X$, and at most $((e/K)n^{2/t-d})^{Kn^d}$ choices of $Y$.
Thus the probability that there is a choice of $X,Y$ such that for every vertex in $X$, all its neighbours are in $Y$,
is at most the product of these, that is
$$\left(Kn^{d-2/t}\right)^{\vare(t-1)n^c} \left((e/\vare)n^{1-c}\right)^{\vare n^c} \left((e/K)n^{2/t-d}\right)^{Kn^d}.$$
The logarithm of this ($L$ say) is 
$$\vare (t-1)n^c\left(\log K +(d-2/t)\log n\right) + \vare n^c \left(\log(e/\vare)+ (1-c) \log n\right) + Kn^d \left(\log (e/K) + (2/t-d)\log n\right).$$
Since $c>d, 1/t$, for sufficiently large $n$ the two terms containing $n^c\log n$ are much larger than the others, and the sum of
their coefficients is $\vare(t-1)(d-2/t) + \vare(1-c)$. This is negative, since $d-1/t < (c-1/t)/(t-1)$; and so
for sufficiently large $n$, $L$ is large and negative, and therefore with high probability, the second bullet of the theorem holds.

Now let $X\subseteq A$ and $Y\subseteq B$, with $|X|=\vare n$ and $|Y|= \lfloor(\vare/e) |B|\rfloor$. 
The probability that for every vertex in $X$, all its neighbours are in $Y$, is at most
$(\vare/e)^{(t-1)\vare n}$.
By \ref{binom} there are at most $(e/\vare)^{\vare n}$ choices of $X$, and at most $2^{n^{2/t}}$ choices of $Y$.
Thus the probability that there is a choice of $X,Y$ such that for every vertex in $X$, all its neighbours are in $Y$,
is at most
$$(\vare/e)^{(t-1)\vare n}(e/\vare)^{\vare n}2^{n^{2/t}}.$$
The logarithm of this ($L$ say) is
$$(t-1)\vare n\log (\vare/e)+ \vare n \log (e/\vare) +n^{2/t}\log 2.$$
The two terms linear in $n$ dominate for large $n$, and the sum of their coefficients is
$$(t-1)\vare \log (\vare/e) + \vare \log (e/\vare)= -(t-2)\vare\log (e/\vare)<0,$$
so
for sufficiently large $n$, $L$ is large and negative, and therefore with high probability, the third bullet of the theorem holds.

Finally, let $v\in B$ and let $X\subseteq A$ with $|X|=n^{1-d}$. The probability that $v$ is adjacent to every vertex in $X$
is at most $\left((t-1)n^{-2/t}\right)^{n^{1-d}}$. The number of choices of $X$ is at most $\left(en^{d}\right)^{n^{1-d}}$, so the probability
that some vertex in $B$ has degree at least $n^{1-d}$ is at most 
$$\left((t-1)n^{-2/t}\right)^{n^{1-d}}\left(en^{d}\right)^{n^{1-d}}n.$$
The logarithm of this is 
$$n^{1-d} \left(\log (t-1) -(2/t) \log n\right) + n^{1-d} (1 + d \log n) + \log n.$$
The $n^{1-d} \log n$ terms dominate, for large $n$, and the sum of their coefficients is $-2/t+ d$; and this is negative since
$d<2/t$. Consequently with high probability, the fourth bullet of the theorem holds.

This proves \ref{bigrandom1}.~\bbox

\begin{thm}\label{bigrandom2}
Let $t\ge 3$ be an integer. Let $0<\vare<1$ be rational, and let $c>d>1/t$, where $d-1/t < (c-1/t)/(t-1)$ and $d<2/t$, 
and $c,d$ are rational.
Let $K>\log (e/\vare)$.
Let $n$ be an integer such that
$n^c, \vare n^c,n^d, n^{1/t}, \vare n, n^{1-d}, (\vare/2) n^d$ are all integers. If $n$ is sufficiently large,
there is a graph with bipartition $B,C$, where $|B|=n^{2/t}$ and $|C|=n$, such that
\begin{itemize}
\item every vertex in $B$ has degree at most $n^{1-d}$;
\item for every $X\subseteq B$ with $|X|\ge Kn^d$, there are more than $(1-\vare) n$ vertices in $C$ that have a neighbour in $X$; and
\item for every $X\subseteq B$ with $|X|\ge (\vare/e) |B|$, there are more than $n-\vare n^c$ vertices in $C$ with a neighbour in $X$.
\item every vertex in $C$ has degree at most $(\vare/2) n^{d}$.
\end{itemize}
\end{thm}
\Proof
Let $B,C$ be disjoint sets of cardinalities $n^{2/t},n$ respectively. For each $v\in B$, choose $n^{1-d}$ vertices in $B$ uniformly
and independent at random (and therefore not necessarily distinct), and add edges to make $v$ adjacent to them. Let $G$
be the graph this constructs. We claim that if $n$ is sufficiently large then with high probability $G$ satisfies the theorem.

Let $X\subseteq B$ with $|X|= Kn^d$, and let $Y\subseteq C$ with $|Y|=(1-\vare)n$. 
The probability that for every vertex in $X$, all its neighbours are in $Y$, is
$$(1-\vare)^{Kn}\le e^{-\vare Kn}.$$
By \ref{binom} there are at most $\left((e/K)n^{2/t-d}\right)^{Kn^d}$ choices of $X$, and at most $(e/\vare)^{\vare n}$ choices of $Y$.
Thus the probability that there is a choice of $X,Y$ such that for every vertex in $X$, all its neighbours are in $Y$, 
is at most 
$$e^{-\vare Kn}\left((e/K)n^{2/t-d}\right)^{Kn^d}(e/\vare)^{\vare n}.$$
The logarithm of this is 
$$-\vare Kn + Kn^d \left(\log (e/K)+(2/t-d) \log n\right) + \vare n \log (e/\vare).$$
For sufficiently large $n$ the terms linear in $n$ dominate, and the sum of their coefficients is 
$-\vare K + \vare \log (e/\vare)$; and this is negative since  $K>\log (e/\vare)$,
so with high probability, the second bullet of the theorem holds.

Now let $X\subseteq B$ and $Y\subseteq C$, with $|X|=\lceil (\vare/e) |B|\rceil$ and $|Y|= n-\vare n^c$. 
The probability that for every vertex in $X$, all its neighbours are in $Y$, is at most
$$\left(1-\vare n^{c-1}\right)^{(\vare/e)n^{2/t}n^{1-d}}=\left(1-\vare n^{c-1}\right)^{(\vare/e)n^{1+2/t-d}} \le e^{-\vare n^{c-1}(\vare/e)n^{1+2/t-d}}=
e^{-(\vare^2/e)n^{c+2/t-d}}.$$
By \ref{binom} there are at most $2^{n^{2/t}}$ choices of $X$, and at most $\left((e/\vare)n^{1-c}\right)^{\vare n^c}$ choices of $Y$.
Thus the probability that there is a choice of $X,Y$ such that for every vertex in $X$, all its neighbours are in $Y$,
is at most
$$e^{-(\vare^2/e)n^{c+2/t-d}}2^{n^{2/t}}\left((e/\vare)n^{1-c}\right)^{\vare n^c}.$$
The logarithm of this is
$$-(\vare^2/e)n^{c+2/t-d} + n^{2/t}\log 2 + \vare n^c \left(\log(e/\vare)+(1-c) \log n\right).$$
Since $c+2/t-d >\max(2/t,c)$, the first term dominates if $n$ is sufficiently large, and so with high probability, 
the third bullet of the theorem holds.

Let $v\in C$ and let $X\subseteq B$ with $|X|=(\vare/2) n^d$. The probability that $v$ is adjacent to every vertex in $X$ is
 $\left(n^{-d}\right)^{(\vare/2) n^d}=n^{-(\vare/2) d n^d}$. 
There are at most $\left((2e/\vare)n^{2/t-d}\right)^{(\vare/2) n^d}$ choices of $X$ by \ref{binom}, and $n$ choices of $v$, so
the probability that there is a choice of $v,X$ such that $v$ is adjacent to every vertex in $X$,
is at most
$$n^{-(\vare/2)  d n^d}\left((2e/\vare)n^{2/t-d}\right){(\vare/2) n^d}n.$$
The logarithm of this is
$$-(\vare/2)  d n^d \log n + (\vare/2) n^d \left(\log (2e/\vare) + (2/t-d)\log n\right) + \log n.$$
The terms in $n^d\log n$ dominate for large $n$, and the sum of their coefficients is 
$$-(\vare/2)  d+  (\vare/2) (2/t-d)=\vare(1/t-d)<0.$$
Consequently with high probability the fourth bullet holds.
This proves \ref{bigrandom2}.~\bbox

\begin{thm}\label{bigrandom3}
Let $0<c\le 1$. If $n$ is sufficiently large, and $\vare n$, $\vare n^c$ are both integers,
there is a graph with bipartition $A,B$, where $|A|=|B|=n$, such that 
every vertex has degree at most $(2/\vare^2)n^{1-c}$, and the blockade
$(A,B)$ is $(\vare n, \vare n^c)$-cohesive.
\end{thm}

We leave the proof to the reader; it is like that of \ref{random}.
The three preceding lemmas are used for the following:

\begin{thm}\label{orderedstar}
Let $t\ge 3$, and let $S_t^+$ be obtained from $S_t$ by ordering its vertex set.
For all $c>1/t$ and all $\vare>0$, there is a graph $G$, and a blockade $(B_1\LL B_{t+1})$ in $G$, 
with local degree less than $\vare W$ and $(\vare W,\vare W^{c})$-cohesive where $W$ is its width, such that there
is no ordered $\mathcal{B}$-rainbow copy of $S_t^+$ in $G$.
\end{thm}

\Proof We may assume that $c,\vare$ are rational, by slightly decreasing them if necessary. 
We call the vertex of $S_t$ of degree $t$ its {\em centre}. From the symmetry we may assume that
the centre is the last in the ordering of $S_t^+$.
Choose $d$ such that $c>d>1/t$, where $d-1/t < (c-1/t)/(t-1)$ and $d<2/t$, and $d$ is rational.
Choose an integer $n$ such that $n^c, n^d, n^{1/t}, \vare n, n^{1-d}, \vare n^d$ are all integers, 
and $n$ is large enough to satisfy each of 
\ref{bigrandom1}, \ref{bigrandom2} and \ref{bigrandom3}.
Choose $K$ as in \ref{bigrandom1} and \ref{bigrandom2}.

Take $t+2$ pairwise disjoint sets $B_0,B_1\LL B_{t+1}$, where $|B_0|=n^{2/t}$ and $B_1\LL B_{t+1}$  all have cardinality $n$.
We attach bipartite graph onto various pairs of the sets $B_0\LL B_{t+1}$ as follows:
\begin{itemize}
\item Let $J_{t+1,0}$ be a copy of the graph of \ref{bigrandom1} with bipartition $B_{t+1}, B_0$. 
\item For $1\le i\le t$, let $J_{0,i}$ be a copy of the graph
of \ref{bigrandom2} with bipartition $(B_0, B_i)$.
\item For $1\le i<j\le t$, let $J_{i,j}$ be a copy of the graph of \ref{bigrandom3} with bipartition $B_i, B_j$.
\end{itemize}
Now for $1\le i<j\le t+1$, and all $u\in B_i$ and $v\in B_j$, add an edge between $u,v$ if they have a common neighbour in $B_0$.
Finally, delete $B_0$; this defines a graph $G$, with a blockade $\mathcal{B}=(B_1\LL B_{t+1})$ of width $W=n$, and we claim it satisfies the 
theorem.

Suppose first that there is an ordered $\mathcal{B}$-rainbow copy of $S_t^+$ in $G$. Thus there exists $v_{t+1}\in B_{t+1}$ 
adjacent in $G$ to some
$v_i\in B_i$ for $1\le i\le t$, such that $v_1\LL v_t$ are pairwise $G$-nonadjacent. From the construction, for $1\le i\le t$ 
every vertex of $B_i$
$G$-adjacent to $v_{t+1}$ is $J_{0,i}$-adjacent to a vertex $w\in B_0$ that is $J_{t+1,0}$-adjacent  to $v_{t+1}$. There are only $t-1$
such vertices, because of the properties of $J_{t+1,0}$; so there exist distinct $i,j\in \{1\LL t\}$ and $w\in B_0$, such that
$v_i, v_j,v_{t+1}$ are all adjacent (in $J_{0,i}, J_{0,j}, J_{t+1,0}$ respectively) to $w$. But then $v_i$ is $G$-adjacent to $v_j$, a contradiction.
This proves that there is no ordered $\mathcal{B}$-rainbow copy of $S_t^+$ in $G$.

To check the local degree of $\mathcal{B}$: each vertex in $B_{t+1}$ is $J_{t+1,0}$-adjacent to at most $t-1$ vertices in $B_0$; 
and each of these neighbours has at most $n^{1-d}$ $J_{0,i}$-neighbours in $B_i$; so each vertex in $B_{t+1}$ has at most 
$(t-1)n^{1-d}$
$G$-neighbours in $B_i$, for $1\le i\le t$, and $(t-1)n^{1-d}<\vare n$ if $n$ is large enough.

Each vertex in $B_i$ is $J_{0,i}$-adjacent to at most $(\vare/2) n^{d}$ vertices in $B_0$, and they have degree less than
$n^{1-d}$ in $J_{t+1,0}$; so each vertex in $B_{i}$ has at most 
$\vare n/2$
$G$-neighbours in $B_{t+1}$, for $1\le i\le t$.

For $1\le i<j\le t$, each vertex in $B_i$ has at most $(2/\vare^2)n^{1-c}$ $J_{i,j}$-neighbours in $B_j$; and in addition,
it is $J_{0,i}$-adjacent to most $(\vare/2) n^{d}$ vertices in $B_0$, and they have degree at most $n^{1-d}$ in $J_{0,j}$.
Consequently each vertex in $B_i$ is $G$-adjacent to at most $(2/\vare^2)n^{1-c} + (\vare/2) n^{d}n^{1-d}$ vertices in $B_j$.
For large $n$, this is less than $\vare n$. The same holds for $j>i$.
Consequently $\mathcal{B}$ has local degree less than $\vare n$.

To check that it is $(\vare n, \vare n^c)$-coherent: first, let $X\subseteq B_{t+1}$ and $Y\subseteq B_i$ where $1\le i\le t$,
with $|X|\ge \vare n$ and $|Y|\ge \vare n^c$. From the properties of $J_{t+1,0}$ there is a set $Z\subseteq B_0$ with
$|Z|\ge (\vare/e) n^{2/t}$, and all its members have
a $J_{t+1,0}$-neighbour in $X$; and from the properties of $J_{0,i}$, more than $n-\vare n^c$ vertices in $B_{i}$
have a $J_{0,i}$-neighbour in $Z$; and consequently some vertex in $Y$ is $G$-adjacent to a vertex in $X$.

Next, let $X\subseteq B_{t+1}$ and $Y\subseteq B_i$ where $1\le i\le t$,
with $|X|\ge \vare n^c$ and $|Y|\ge \vare n$. By a similar argument it follows that $X,Y$ are not anticomplete in $G$.

Finally, let $1\le i<j\le t$, and let $X\subseteq B_{i}$ and $Y\subseteq B_j$, with $|X|\ge \vare n$ and $|Y|\ge \vare n^c$.
From the properties of $J_{i,j}$, $X,Y$ are not anticomplete in $J_{i,j}$ and hence not in $G$. This proves \ref{orderedstar}.~\bbox

\end{document}